\documentclass[12pt]{amsart}
\usepackage{latexsym,amsmath,amsfonts,amscd,amssymb}
\newcommand{\n}{\noindent}
\newtheorem{theorem}{\bf Theorem}[section]
\newtheorem{lemma}[theorem]{\bf Lemma}
\newtheorem{remark}[theorem]{\bf Remark}

\newtheorem{proposition}[theorem]{\bf Proposition}
\newtheorem{corollary}[theorem]{\bf Corollary}
\newtheorem{definition}[theorem]{\bf Definition}

\newtheorem{cond}[theorem]{\bf Conditions}
\newcommand{\be}{\begin{equation}}
\newcommand{\ee}{\end{equation}}

\newfont{\bfc}{cmbsy10 scaled 1200}  % bold face calligraphic
\newfont{\dr}{msbm10 scaled \magstep1}  %letra doble raya
\newfont{\sdr}{msbm8}  % small letra doble raya
\newfont{\gl}{eufm10 scaled \magstep1}  % german letters as well

\DeclareFontFamily{OT1}{rsfs}{}
\DeclareFontShape{OT1}{rsfs}{n}{it}{<->rsfs10}{}
\DeclareMathAlphabet{\curly}{OT1}{rsfs}{n}{it}
\renewcommand{\a}{\alpha}
\renewcommand{\b}{\beta}
\newcommand{\p}{\phi}
\renewcommand{\l}{\lambda}

 \newcommand{\CC}{{\mathbb C}}

 \newcommand{\HH}{{\mathbb H}}
 
 \newcommand{\PP}{{\mathbb P}}
 
 \newcommand{\RR}{{\mathbb R}}

 \newcommand{\VV}{{\mathbb V}}
 
 \newcommand{\ZZ}{{\mathbb Z}}

 \newcommand{\cF}{{\mathcal F}}

 \newcommand{\cO}{{\mathcal O}}

\newcommand{\qu}{/\kern-.7ex/}
\newcommand{\exh}{\to\kern-1.8ex\to}
 \newcommand{\ra}{\rightarrow}
 \newcommand{\lra}{\longrightarrow}

\newcommand{\Aut}{\operatorname{Aut}}
\newcommand{\End}{\operatorname{End}}
\newcommand{\GL}{\operatorname{GL}}
\newcommand{\GCD}{\operatorname{GCD}}
\newcommand{\Gr}{\operatorname{Gr}}
\newcommand{\Hom}{\operatorname{Hom}}
\newcommand{\Ext}{\operatorname{Ext}}

\newcommand{\Ker}{\operatorname{Ker}}

\renewcommand{\Im}{\operatorname{Im}}
\newcommand{\Pic}{\operatorname{Pic}}

\newcommand{\rk}{\operatorname{rk}}

\newcommand{\codim}{\operatorname{codim}}

\newcommand{\gr}{\operatorname{gr}}

\newcommand{\ndk}{(n,d,k)}
\newcommand{\x}{\times}
\newcommand{\ox}{\otimes}
\renewcommand{\tilde}{\widetilde}
\begin{document}
\title{Coherent systems and Brill-Noether theory}
\thanks{All authors are members of the research group VBAC
(Vector Bundles on Algebraic Curves), which is partially supported by
EAGER (EC FP5 Contract no.~HPRN-CT-2000-00099) and by EDGE (EC FP5
Contract no. HPRN-CT-2000-00101). Support was also received from the
Acciones Integradas Programme (HB 1998-0006). The first author was
partially supported by the National Science Foundation under grant
DMS-0072073. The fourth author would like to thank the Isaac Newton
Institute, Cambridge and the organisers of the HDG programme for their
hospitality during the completion of work on this paper.}

\subjclass{14D20, 14H51, 14H60}
\date{29 May 2002}

\keywords{Algebraic curves, moduli of vector bundles, coherent
systems, Brill-Noether loci}

\author{S.B. Bradlow}
  \address{Department of Mathematics\\University of Illinois\\
  Urbana\\IL 61801\\USA}
  \email{bradlow@math.uiuc.edu}
\author{O. Garc\'{\i}a-Prada}
   \address{Departamento de Matem{\'a}ticas \\
   Facultad de Ciencias \\ Universidad Aut{\'o}noma de Madrid \\
   28049 Madrid \\ Spain}
   \email{oscar.garcia-prada@uam.es}
\author{V. Mu\~noz}
  \address{Departamento de Matem\'aticas \\
  Facultad de Ciencias \\ Universidad Aut\'onoma de Madrid
  \\ 28049 Madrid \\ Spain}
  \email{vicente.munoz@uam.es}
\author{P.E. Newstead}
   \address{Department of Mathematical Sciences \\
   University of Liverpool \\ Peach Street \\
   Liverpool L69 7ZL \\ UK}
   \email{newstead@liverpool.ac.uk}

%\begin{document}
\maketitle
%\begin{abstract}
%\end{abstract}
%\tableofcontents

\hspace{208pt}\begin{minipage}[t]{220pt}\scriptsize Suplico a vuesa merced,
se\~nor don Quijote, que mire bien y especule con cien ojos lo que hay
all\'a dentro: quiz\'a habr\'a cosas que las ponga yo en el libro de mis
{\em Transformaciones} (El ingenioso hidalgo don Quijote de la Mancha,
Book 2, Chapter XXII)

I beg you, don Quixote sir: look carefully, inspect with a hundred eyes what
you see down there. Who knows, maybe you will find something that I can put in
my book on {\em Transformations}.
\end{minipage}
%%%%%%%%%%%%%%%%%%%%%%%
\section{Introduction}
\label{sec:intro}
%%%%%%%%%%%%%%%%%%%%%%%

\renewcommand{\theenumi}{\roman{enumi}}
\renewcommand{\labelenumi}{(\theenumi)}
Augmented algebraic vector bundles often have moduli spaces which
depend not only on the topological type of the augmented bundle,
but also on an additional parameter. The result is that the moduli
spaces occur in discrete families. First exploited by Thaddeus in
a proof of the Verlinde formula, this phenomenon has been
responsible for several interesting applications (cf.~\cite{BeDW},
\cite{BGG}, \cite{BG2}). In this paper we examine the augmented
bundles known as coherent systems and discuss the use of their
moduli spaces as a tool in Brill-Noether theory.

By a {\it coherent system} on an algebraic variety (or scheme) we
mean an algebraic vector bundle together with a linear subspace of
prescribed dimension of its space of sections. As such it is an
example of an augmented bundle. Introduced in \cite{KN}, \cite{RV}
and \cite{LeP}, there is a notion of stability which permits the
construction of moduli spaces. This notion depends on a real
parameter, and thus leads to a family of moduli spaces.

That these moduli spaces are related to Brill-Noether loci follows
almost immediately from the definitions.  The Brill-Noether loci are
natural subvarieties within the moduli spaces of stable bundles over
an algebraic curve, defined by the condition that the bundles have at
least a prescribed number of linearly independent sections. A similar
condition defines Brill-Noether loci in the moduli spaces of
(S-equivalence classes of) semistable bundles. But any bundle which
occurs as part of a coherent system must evidently have at least a
prescribed number of linearly independent sections. Conversely,  a
bundle with a prescribed number of linearly independent sections
determines, in a natural way, a coherent system.

In order to convert these observations into a precise relationship
between the coherent systems moduli spaces and the Brill-Noether
loci, one extra tool is needed, namely a precise relationship between
bundle stability and coherent system stability. In general, for an
arbitrary choice of the coherent system stability parameter, no
such relationship exists. However, for values of the parameter close
to $0$, the required relationship holds and there is a map from the
coherent systems moduli space to the semistable Brill-Noether
locus. While this map is not necessarily surjective, it does
include the entire stable Brill-Noether locus in its image.  It
is via this map that information about the coherent systems moduli
spaces can be applied to answer questions about higher rank
Brill-Noether theory.

In \cite{BG2} the first two authors initiated a programme to do
just this, i.e.~to use coherent systems moduli spaces to study
higher rank Brill-Noether theory. There the goals were limited to
explaining results of \cite{BGN} from the perspective of coherent
systems. This turned out to require only a limited understanding
of the coherent systems moduli spaces. In this paper we build on
the foundation laid in \cite{BG2}.

We study the moduli spaces for coherent systems $(E,V)$\
consisting of an algebraic vector bundle $E$ together with a
linear subspace $V$ of its space of
sections. While not required by the definitions, we consider only
the case of bundles over a smooth irreducible projective algebraic
curve $X$ of genus $g$.
The {\it type} of the coherent system is defined by a triple of integers
$(n,d,k)$ giving the rank of $E$, the degree of $E$, and the
dimension of the subspace $V$.

The infinitesimal study of coherent systems follows a standard pattern and is
summarised in section \ref{sec:infinitesimal}. This allows us in particular
to identify the Zariski tangent space to each moduli space at any point,
and to show that every irreducible component of every moduli space has dimension
at least equal to a certain number $\b(n,d,k)$, called the {\em
Brill-Noether number} and often referred to as the {\em expected dimension}.

For each type $(n,d,k)$, there is a family of
moduli spaces. While each such family of moduli spaces has some properties which
depend on
$\ndk$, there are some features that are common to all types. In
particular:
\begin{itemize}
\item {\it The families have only a finite number of distinct members.}
The different members in the family correspond to different values
for the real parameter $\a$ in the definition of stability. As
$\a$ varies, the stability condition changes only as $\a$ passes
through one of a discrete set of points in the real line. In some
cases (if $k<n$) the range for the parameter is a finite interval,
in which case it follows automatically that the family has only
finitely many distinct members. However, even in the cases for
which the range of the parameter is infinite, it turns out that
there can be only a finite number of distinct moduli spaces. This
is a consequence of the {\it stabilisation theorem} Proposition
\ref{prop:flips}.
\item {\it  The families are ordered, and the coherent systems in
the terminal member are as simple as possible.} The ordering comes
from the fact that the moduli spaces are labelled by intervals on the real
line. By the terminal member of the family we mean the moduli
space corresponding to the last of these intervals (in the natural ordering
of $\RR$).  In
section \ref{sec:alpha-large} we analyse the coherent systems corresponding
to
the
points in this terminal moduli space. While the specifics depend
on the type of the coherent system, in all cases we find that the
structure of these coherent systems is (in a suitable sense) the
best possible.
\end{itemize}

The most obvious type-dependent feature is the description of the
terminal moduli space. This divides naturally into distinct cases,
according to whether $k<n$, $k=n$ or $k> n$.  The case $k<n$ was
discussed in some detail in \cite{BG2}. The results are summarized in
section \ref{subsec:alpha-large.k<n}. When $k\ge n$, we show (in
section \ref{subsec:quot}) how to relate the terminal moduli space to
a Grothendieck Quot scheme of quotients of the trivial bundle of rank
$k$. Denoting the terminal moduli space of stable (respectively
semistable) objects by $G_L$ (respectively $\tilde{G_L}$), we prove

\begin{theorem}\label{thm:main1} {\bf [Theorem \ref{G_L(k=n)}]}
Let $k=n\ge2$.  If $\tilde{G}_L$ is non-empty then $d=0$ or $d\geq
n$. For $d>n$,
  $\tilde{G}_L$ is irreducible and $G_L$ is smooth  of the expected
dimension $dn-n^2+1$.
  For $d=n$, $G_L$ is empty and $\tilde{G}_L$ is irreducible
of dimension $n$ (not of the expected dimension). For $d=0$,
$G_L$ is empty and $\tilde{G}_L$ consists of a single point.
\end{theorem}

The non-emptiness of $G_L$ for the case $k>n$ is not so obvious
and is related to the non-emptiness of Quot schemes. However, in
this case there is a duality construction that relates coherent
systems of types $(n,d,k)$ and $(k-n,d,k)$. If the parameter is
large these ``dual'' moduli spaces are birationally equivalent. A
similar idea has been considered for small $\a$ by Butler
\cite{Bu2}, but the construction seems to be more natural  for
large $\a$ and turns out to be  an important tool to prove
non-emptiness for large values of the parameter. For instance if
$k=n+1$, the  non-emptiness is given by the classical rank 1
Brill-Noether theory, i.e. we get

\begin{theorem} \label{thm:main2} {\bf [Theorem
\ref{thm:dual-span}]}
Suppose that the curve $X$ is generic and that $k=n+1$. Then $G_L$
is non-empty if and only if $\b=g-(n+1)(n-d+g)\ge0$. Moreover
$G_L$ has dimension $\b$ and it is irreducible whenever $\b>0$.
\end{theorem}

In addition to these `absolute' results about the terminal moduli
spaces, we also give `relative' results which
characterize the differences between moduli spaces within a given
family. We compare the moduli spaces and identify subvarieties
within which the differences are localised. In section
\ref{sec:critical} we give some general results estimating the
codimensions of these subvarieties.

With a view to applications, we examine a number of special cases
in which either $k$ or $n$ (or both) are small. In these cases,
discussed in sections \ref{sec:k=1} - \ref{sec:k=3}, we can get more
detailed results,
especially for the codimension estimates on the difference loci
between moduli spaces within a family.

Having amassed all this information about the coherent systems
moduli spaces, we end with some applications to Brill-Noether
theory in section \ref{sec:applications}. In all cases the strategy is the same:
starting with
information about $G_L$, and using our results about the relation
between different moduli spaces within a given family, we deduce
properties of the moduli space corresponding to the smallest values
of the stability parameter. This is then passed down to the
Brill-Noether loci using the morphism from the coherent systems
moduli space to the moduli space of semistable bundles. This allows us,
for example,
\begin{itemize}
\item {\bf [Theorems \ref{thm:bn2} and \ref {thm:bn3}]} to prove
the irreducibility
of
the Brill-Noether loci for $k=1, 2, 3$, and
\item {\bf [Theorem \ref{thm:bn4}]} to compute the Picard group of the
smooth part of
the
Brill-Noether locus for $k=1$.
\end{itemize}

While the irreducibility result was previously known for
$k=1$ and any $d$, and for $k=2,3$, $k<n$ and $d<\min\{2n,n+g\}$,
our theorems have no such restrictions. These results should be
regarded as a sample of what can be done. Our methods are certainly
applicable more widely and we propose to pursue this in future
papers.

Throughout the paper $X$ denotes a fixed smooth irreducible projective
algebraic curve of genus $g$ defined over the complex numbers. Unless
otherwise stated,
we make no assumption about  $g$. For simplicity
we shall write $\cO$ for $\cO_X$ and $H^0(E)$ for $H^0(X,E)$. We shall
consistently denote the ranks of bundles $E,E',E_1\ldots$ by
$n,n',n_1\ldots$,
their degrees by $d,d',d_1\ldots$ and the dimensions of subspaces
$V,V',V_1\ldots$
of their spaces of sections by $k,k',k_1\ldots$.

%%%%%%%%%%%%%%%%%%%%%%%%%%%%%%%%%%%%%%%%%%%%%%%%%%
\section{Definitions and basic facts}
\label{sec:definitions}
%%%%%%%%%%%%%%%%%%%%%%%%%%%%%%%%%%%%%%%%%%%%%%%%%%

%%%%%%%%%%%%%%%%%%%%%%%%%%%%%%%%%%%%%%%%%%%%%%%%%%%%
\subsection{Coherent systems and their moduli spaces}
\label{subsec:coh}
%%%%%%%%%%%%%%%%%%%%%%%%%%%%%%%%%%%%%%%%%%%%%%%%%%%%

Recall (cf. \cite{LeP}, \cite{KN}) that a coherent system $( E,
V)$\ on  $X$\ of type $(n,d,k)$ consists of an algebraic vector
bundle $ E$ over $X$\ of rank $n$\ and degree $d$, and  a linear
subspace $V\subseteq H^0(E)$\ of dimension $k$.  Strictly
speaking, it is better to consider triples $( E,\VV,\p)$\
where $\VV$\ is a dimension $k$\ vector space and $\p:\VV
\ox\cO\ra E$\ is a sheaf map such that the induced map
$H^0(\p):\VV\ra H^0(E)$\ is injective. The linear space $V\subseteq
H^0(E)$\ is then the image $H^0(\p)(\VV)$. Under the natural concepts
of isomorphism, isomorphism classes of such triples are in bijective
correspondence with isomorphism classes of coherent systems. We will
usually use the simpler notation $( E, V)$, but occasionally it is
helpful to use the longer one. For a
summary of basic results about coherent systems (and other related
augmented bundles) we refer the reader to \cite{BDGW}.

By introducing a suitable definition of stability, one can
construct moduli spaces of coherent systems. The correct notion
(i.e. the one dictated by Geometric Invariant Theory) depends on a
real parameter $\a$, which a posteriori must be non-negative (cf.
\cite{KN}). In the situation under consideration (i.e. where $E$\
is a vector bundle over a smooth algebraic curve), the definition
may be given as follows.
\begin{definition} \label{def:stable}
\begin{em}
Fix $\a\in\RR$. Let $( E,V)$\ be a coherent system of type $\ndk$.
The $\a$-{\em slope} $\mu_{\a}(E,V)$ is defined by
 $$
 \mu_{\a}( E,V)=\frac{d}{n}+\a\frac{k}{n}\ .
 $$
We say $( E,V)$\ is $\a$-{\em stable} if
 $$
 \mu_{\a}(E',V')<\mu_{\a}(E,V)
 $$
for all proper subsystems $(E',V')$ \textup{(}i.e. for every
non-zero subbundle $E'$ of $E$ and every subspace $V'\subseteq
V\cap H^0(E')$ with $(E',V')\ne(E,V)$\textup{)}. We define
$\a$-{\em semistability} by replacing the above strict inequality with a
weak inequality. A coherent system is called $\a$-{\em polystable} if it
is the direct sum of $\a$-stable coherent systems of the same
$\a$-slope.\end{em}
\end{definition}

Sometimes it is necessary to consider a larger class of objects
than coherent systems in which one replaces $E$ by a general
coherent sheaf and $H^0(\p):\VV\ra H^0(E)$  is not necessarily
injective. By doing so one obtains an abelian category \cite{KN}.
One can easily extend the definition of $\a$-stability to this
category. It turns out, however, that $\a$-semistability forces
$E$ to be locally free and $H^0(\p)$ to be injective, and hence
 $\a$-semistable objects in this category can be identified with
$\a$-semistable  coherent systems up to an appropriate definition of
isomorphism.
One has the following result.

\begin{proposition}\label{prop:filtration}(\cite[Corollary 2.5.1]{KN})
The $\a$-semistable coherent systems of any fixed $\a$-slope form a
noetherian and artinian abelian category in which the simple objects
are precisely the $\a$-stable systems. In particular the following
statements
hold.

\begin{enumerate}
\item {\em (Jordan-H\"older Theorem)} For any $\a$-semistable coherent
system $(E,V)$, there exists a   filtration
 by $\a$-semistable coherent systems $(E_j,V_j)$,
 $$
 0=(E_0,V_0)\subset (E_1,V_1)\subset ...\subset (E_m,V_m)=(E,V),
 $$
 with $(E_j,V_j)/(E_{j-1},V_{j-1})$ an $\a$-stable coherent system and
 $$
 \mu_\a((E_j,V_j)/(E_{j-1},V_{j-1}))=\mu_\a(E,V)
 \;\;\; \mbox{for}\;\; 1\leq j\leq m.
 $$
\item If $(E,V)$ is an $\a$-stable coherent system, then
$\End(E,V)\cong\CC$.
\end{enumerate}
\end{proposition}

 Any filtration as in (i) is called
a {\em Jordan-H\"older} filtration of $(E,V)$. It is not
necessarily unique, but the associated graded object is uniquely
determined by $(E,V)$.

\begin{definition}\label{def:polystable}\begin{em}
We define the {\em graduation} of $(E,V)$ to be the $\a$-polystable
coherent system
 $$
 \gr(E,V)=\bigoplus_j (E_j,V_j)/(E_{j-1},V_{j-1}).
 $$
Two $\a$-semistable coherent systems $(E,V)$ and $(E',V')$ are
said to be $S$-equivalent if $\gr(E,V)\cong\gr(E',V')$.\end{em}
\end{definition}
We shall denote the moduli space of   $\a$-stable coherent systems
of type $(n,d,k)$ by $G(\a)=G(\a;n,d,k)$, and the moduli space
of S-equivalence classes of  $\a$-semistable coherent systems
of type $(n,d,k)$
 by $\tilde{G}(\a)= \tilde{G}(\a;n,d,k)$. The
moduli space $\tilde{G}(\a)$ is a projective variety which
contains $G(\a)$ as an open set.

Now suppose that $k\geq 1$. By applying the $\a$-semistability
condition for $(E,V)$ to the subsystem $(E,0)$ one obtains that
$\a\geq 0$. This means that there are no semistable coherent
systems for negative values of $\a$. For $\a=0$, $(E,V)$ is
$0$-semistable if and only if $E$ is semistable. For $k\geq 1$
there are no 0-stable coherent systems.
\begin{definition}\label{def:cv}\begin{em}
We say that  $\a>0$\  is a  {\em virtual critical value} if
 it is numerically possible to have a proper subsystem
$(E',V')$\ such that $\frac{k'}{n'}\ne\frac{k}{n}$\ but $\mu_{\a}(
E',V')=\mu_{\a}( E,V)$. We also regard $0$ as a virtual critical
value. If there is a  coherent system $(E,V)$ and a subsystem
$(E',V')$ such that this actually holds, we say that $\a$ is an
{\em actual critical value}.\end{em}
\end{definition}
It follows from this  (cf. \cite{BDGW}) that, for coherent systems
of type $(n,d,k)$, the non-zero virtual critical values of $\a$\ all lie in
the set
 $$
 \{\  \frac{nd'-n'd}{n'k-nk'}\ |\  0\le k'\le k\ ,\
 0<n'< n\ ,\ n'k\ne nk'\ \}\cap (0,\infty)\ .
 $$
We say that $\a$ is {\em generic} if it is not critical. Note
that, if $\GCD(n,d,k)=1$ and $\a$ is generic, then
$\a$-semistability is equivalent to $\a$-stability. If we label
the critical values of $\a$\ by $\a_i$, starting with $\a_0=0$, we
get a partition of the $\a$-range into a set of intervals
$(\a_{i},\a_{i+1})$. For numerical reasons it is clear that within
the interval $(\a_{i},\a_{i+1})$ the property of $\a$-stability is
independent of $\a$, that is if $\a,\a'\in (\a_{i},\a_{i+1})$,
$G(\a)= G(\a')$. We shall denote this moduli space by
$G_i=G_i(n,d,k)$. The construction of moduli spaces thus yields
one moduli space $G_i$ for the interval $(\a_i,\a_{i+1})$. If
$\GCD(n,d,k)\ne1$, one can define similarly the moduli spaces
$\tilde{G}_i$ of semistable coherent systems. The GIT construction
of these moduli spaces has been given in \cite{LeP} and \cite{KN}.
A previous construction for $G_0$ had been given in \cite{RV} and
in \cite{Be} for big degrees. When $k=1$ the moduli space of
coherent systems is equivalent to the moduli space of vortex pairs
studied in \cite{B,BD1,BD2,G,HL1,HL2,Th}.

The relationship between the semistability of a coherent system and the
underlying vector bundle is given by the following (cf.
\cite{BDGW}, \cite{KN}).

\begin{proposition} \label{small-alpha} Let $\a_1$ be the first critical
value
after $0$ and let $0<\a<\a_1$. Then

\begin{enumerate}
\item $(E,V)$ $\a$-stable implies $E$ semistable;
\item $E$ stable implies $(E,V)$
$\a$-stable.
\end{enumerate}
\end{proposition}

%%%%%%%%%%%%%%%%%%%%%%%%%%%%%%%%%%%%%%%%%%%%%%%%%%%%%%%%%%%%
\subsection{ Brill-Noether loci}\label{subsec:basicBN}
%%%%%%%%%%%%%%%%%%%%%%%%%%%%%%%%%%%%%%%%%%%%%%%%%%%%%%%%%%%%

\begin{definition} \label{def:BN}\begin{em}
Let $X$ be an algebraic curve, and let $M(n,d)$
be the moduli space of stable bundles of rank $n$ and degree $d$. Let
$k\geq 0$. The {\em Brill-Noether loci} of stable bundles  are defined by
 $$
 B(n,d,k):=\{E\in M(n,d)\; |\; \dim H^0(E)\geq k\}.
 $$\end{em}
\end{definition}
Similarly one defines the  Brill-Noether loci of semistable bundles
 $$
 \tilde{B}(n,d,k):=\{[E]\in \tilde{M}(n,d)\; |\; \dim H^0(\gr(E))\geq k\},
 $$
where $\tilde{M}(n,d)$ is the moduli space of S-equivalence classes
of semistable bundles, $[E]$ is the S-equivalence class of $E$
and $\gr(E)$ is the polystable bundle
defined by a Jordan-H\"older filtration of $E$.

The spaces $B(n,d,k)$ and $\tilde{B}(n,d,k)$ have previously been denoted by
$W_{n,d}^{k-1}$ and $\tilde{W}_{n,d}^{k-1}$, but we have chosen to change
the classical notation in an attempt to get rid of the $k-1$, which
in the arbitrary rank case does not make much sense. (In fact the same loci
have also been denoted by $W_{n,d}^k$ and $\tilde{W}_{n,d}^k$, but, while
logical,
this seems a little confusing!)

By semicontinuity, the Brill-Noether loci are closed subschemes of
the appropriate moduli spaces. The main object of Brill-Noether
theory is the study of these subschemes, in particular questions
related to their non-emptiness, connectedness, irreducibility,
dimension, and topological and geometric structure. It is in
particular not difficult to describe them as determinantal loci,
from which one obtains the following general result. We begin with
a definition.

\begin{definition}\label{def:bnnumber}\begin{em}
For any $(n,d,k)$, the {\em Brill-Noether number} $\b(n,d,k)$ is defined by
\begin{equation}\label{bnnumber}
\b(n,d,k)=n^2(g-1)+1-k(k-d+n(g-1)).
\end{equation}\end{em}
\end{definition}

\begin{theorem}\label{thm:BN}
If $B(n,d,k)$ is non-empty and $B(n,d,k)\neq M(n,d)$, then
\begin{itemize}
\item
every irreducible component $B$ of $B(n,d,k)$ has dimension
 $$
 \dim B\geq \b(n,d,k),
 $$
\item $B(n,d,k+1)\subset Sing B(n,d,k)$,
\item the tangent space of  $B(n,d,k)$  at a point $E$ with $\dim
H^0(E)=k$ can be identified with the dual of the cokernel of
the {\em Petri map}
\begin{equation}\label{petri}
H^0(E)\ox H^0(E^*\ox K)\lra H^0(\End E \ox K)
\end{equation}
(given by multiplication of sections),
\item
$B(n,d,k)$ is smooth of dimension $\b(n,d,k)$ at $E$ if and only if
the Petri map is injective.
\end{itemize}
\end{theorem}

 For details, see for example \cite{M2}.

When $n=1$, $M(n,d)$ is just $J^d$, the Jacobian of $X$ consisting of
line bundles of degree $d$,
and the Brill-Noether loci are the classical ones for which a thorough
modern presentation is given in \cite{ACGH}. In particular we have the
following results.
\begin{itemize}
\item
If $\b(1,d,k)\ge0$, then $B(1,d,k)$ is non-empty.
\item
If $\b(1,d,k)>0$, then $B(1,d,k)$ is connected.
\item
For a generic curve $X$ and $n=1$, the Petri map is always injective. Hence
\begin{itemize}
\item $B(1,d,k)$ is smooth outside $B(1,d,k+1)$.
\item
$B(1,d,k)$ has dimension $\b(1,d,k)$ whenever it is non-empty and
not equal to $M(n,d)$.
\item
$B(1,d,k)$ is irreducible if $\b(1,d,k)>0$.
\end{itemize}
\end{itemize}

None of these statements is true for $n\ge2$ (see, for example,
\cite{T1,T2,BGN,BeF,Mu}).

Rather than referring repeatedly to a {\em generic} curve, we prefer to
use the following more precise term.

\begin{definition}\label{def:Petri}\begin{em}
A curve $X$ is called a {\em Petri curve} if the Petri map
 $$
 H^0(L)\ox H^0(L^*\ox K)\lra H^0(K)
 $$
is injective for every line bundle $L$ over $X$.\end{em}
\end{definition}

One may note that any curve of genus $g\le2$ is Petri, the simplest
examples of non-Petri curves being hyperelliptic curves with $g\ge3$.
There is currently no sensible generalisation of Definition \ref{def:Petri}
to higher rank. Indeed, at least for $g\ge6$, there exist stable bundles
$E$ on Petri curves for which the Petri map (\ref{petri})
is not injective (see \cite[\S5]{T1}).
Moreover the condition of Definition \ref{def:Petri} is not sufficient to
determine even the non-emptiness of Brill-Noether loci in higher rank (see
\cite{M2,Mu,V}).

%%%%%%%%%%%%%%%%%%%%%%%%%%%%%%%%%%%%%%%%%%%%%%%%%%%%%%%%%%%%%%%%%%%%%%%%%%%
\subsection{Relationship between $B(n,d,k)$ and
$G_0$}\label{subsec:alpha-small}
%%%%%%%%%%%%%%%%%%%%%%%%%%%%%%%%%%%%%%%%%%%%%%%%%%%%%%%%%%%%%%%%%%%%%%%%%%%%

The relevance of the moduli spaces of coherent systems in relation
to Brill-Noether theory is given by Proposition
\ref{small-alpha}. The assignment  $(E,V)\mapsto E$ defines a map
\begin{equation}\label{bnmap}
G_0(n,d,k)\lra \tilde{B}(n,d,k),
\end{equation}
which is one-to-one over $B(n,d,k)-B(n,d,k+1)$ and whose image contains
$B(n,d,k)$. When $\GCD(n,d,k)\ne1$, this map can be extended to
\begin{equation}\label{bnmap2}
\tilde{G}_0(n,d,k)\lra \tilde{B}(n,d,k).
\end{equation}

Even (\ref{bnmap2}) may fail to be surjective. This happens
for example (as observed in \cite{BG2}) when $d=0$, $0<k<n$. In this case
$\tilde{G}_0=\emptyset$ but $\tilde{B}$ is non-empty. On the other hand, if
$(n,d)=1$, the loci $B$ and $\tilde{B}$ coincide and (\ref{bnmap}) is
surjective.

Even when $(n,d)\ne1$, we may be able to obtain information about $B$ from
properties of $G_0$. For example, if $G_0(n,d,k)$ is non-empty, then
certainly $\tilde{B}(n,d,k)$ is non-empty. Moreover, if $B(n,d,k)$ is
non-empty and $G_0(n,d,k)$ is irreducible, then $B(n,d,k)$ is also
irreducible.

We are therefore interested in studying $G_0=G_0(n,d,k)$. Our approach to
this consists of having
\begin{itemize}
\item
a detailed description of at least one (usually large $\a$) moduli space,
\item
a thorough understanding of the ``flips'' to go from $G_i$ to $G_{i-1}$,
until we get to
$G_0$. The meaning of `thorough' can vary,
depending on the application. For instance, for non-emptiness
questions all we require are the codimensions of the flip loci, or at
least sufficiently good estimates thereof.
\end{itemize}

In the case $n=1$, everything is much simpler.
The concept of stability is vacuous and independent of $\a
\in (0,\infty)$. We shall therefore denote the moduli space of
coherent systems by $G(1,d,k)=G(\a;1,d,k)$. It consists of
coherent systems $(L,V)$ such that $L$ is a line bundle of
degree $d$ and $V\subset H^0(L)$ is any subspace of dimension $k$.
These spaces have been studied classically (see for example \cite{ACGH},
where $G(1,d,k)$ is denoted by $G_d^{k-1}$). The map (\ref{bnmap}) becomes
\begin{equation}\label{bnmap3}
G(1,d,k)\lra B(1,d,k)
\end{equation}
and is always surjective. The fibre of (\ref{bnmap3}) over $L$ can be
identified with the Grassmannian $\Gr(k,h^0(L))$.

When $X$ is a Petri curve, we have
\begin{itemize}
\item $G(1,d,k)$ is non-empty if and only if $\b=g-k(k-d+g-1)\ge0$,
\item If $\b\ge0$, $G(1,d,k)$ is smooth of dimension $\b$,
\item If $\b>0$, $G(1,d,k)$ is irreducible.
\end{itemize}

Comparing this with section \ref{subsec:basicBN}, we see that, for
$X$ a Petri curve, $G(1,d,k)$ provides a desingularisation of
$B(1,d,k)$ whenever $B(1,d,k)\ne J^d$. In higher rank, if
$\GCD(n,d,k)=1$, $\b(n,d,k)\le n^2(g-1)$ and $G_0(n,d,k)$ is smooth
and irreducible, and
$B(n,d,k)$ is non-empty, then the map (\ref{bnmap}) is a
desingularisation of the closure of $B(n,d,k)$. We shall see that,
for $k=1$, all these conditions hold (see section
\ref{sec:applications}). Note that it was proved in \cite{RV} that
$G_0(n,n(g-1),1)$ is a desingularisation of
$\tilde{B}(n,n(g-1),1)$, which coincides with the generalised
theta-divisor in $M(n,n(g-1))$.

%%%%%%%%%%%%%%%%%%%%%%%%%%%%%%%%%%%%%%%%%%%%%%
\section{Infinitesimal study and extensions}
\label{sec:infinitesimal}
%%%%%%%%%%%%%%%%%%%%%%%%%%%%%%%%%%%%%%%%%%%%%%%%

The infinitesimal study of the moduli space of coherent systems as
well as the study of extensions of coherent systems is carried out in
\cite{He,LeP} (see also \cite{Th,RV}). We review here the main
results and refer to these papers, in particular for omitted proofs.

Given two coherent systems $(E,V)$ and $(E',V')$ one
defines the groups
 $$
 \Ext^q((E',V'),(E,V)),
 $$
and considers the long exact sequence (\cite[Corollaire 1.6]{He})
\begin{equation}\label{long-exact}
\begin{array}{ccccccc}
 0& \lra & \Hom((E',V'),(E,V))   & \lra &  \Hom(E',E) & \lra & \Hom(V',
H^0(E)/V) \\
 & \lra & \Ext^1((E',V'),(E,V)) & \lra &  \Ext^1(E',E) & \lra &  \Hom(V',
H^1(E))\\
 & \lra & \Ext^2((E',V'),(E,V)) & \lra &   0.           &      &
\end{array}
\end{equation}

Notice that since we are on a curve $\Ext^2(E',E)=0$. Also,
since $E'$ is a vector bundle,
 $$
 \Ext^1(E',E)=H^1({E'}^*\ox E).
 $$

We can now apply this to the study of infinitesimal deformations of
the moduli space of coherent systems as well as to the study of
extensions of coherent systems.

%%%%%%%%%%%%%%%%%%%%%%%
\subsection{Extensions}\label{subsec:extensions}
%%%%%%%%%%%%%%%%%%%%%%%

We will have to deal later with extensions of coherent systems
arising from the one-step Jordan-H\"older filtration of a
semistable coherent system. By standard results on abelian
categories, we have

\begin{proposition} \label{prop:extensions}
Let $(E_1,V_1)$ and
$(E_2,V_2)$ be two coherent systems on $X$. The space of
equivalence classes of extensions
 $$
 0\lra (E_1,V_1)\lra (E,V) \lra (E_2,V_2)\lra 0
 $$
is isomorphic to $\Ext^1((E_2,V_2),(E_1,V_1))$. Hence the quotient
of the space of non-trivial extensions by the natural action of
$\CC^*$ can be identified with the projective space
$\PP(\Ext^1((E_2,V_2),(E_1,V_1)))$.
\end{proposition}

\begin{proposition}\label{prop:C21}
Let $(E_1,V_1)$ and $(E_2,V_2)$ be two coherent systems
on $X$ of types $ (n_1,d_1,k_1)$ and $(n_2,d_2,k_2)$ respectively.
Let  ${\HH}_{21}^0:=\Hom((E_2,V_2),(E_1,V_1))$ and
 $\HH^2_{21}:=\Ext^2((E_2,V_2),(E_1,V_1))$. Then
\begin{equation}\label{C21}
\dim \Ext^1((E_2,V_2),(E_1,V_1))= C_{21}+\dim {\HH}^0_{21}+ \dim
{\HH}^2_{21},
\end{equation} where
\begin{eqnarray}
C_{21} &:= & k_2\chi (E_1)-\chi(E_2^\ast\ox E_1)- k_1k_2
\nonumber\\
       & = &  n_1n_2(g-1)-d_1n_2+d_2n_1+k_2d_1-k_2n_1(g-1)-k_1k_2.
                      \label{dim-ext}
\end{eqnarray}
Moreover,
   \begin{equation} \label{obs-ext}
     {\HH}_{21}^2=\Ker(H^0(E_1^*\ox K)\ox V_2 \to
     H^0(E_1^* \ox E_2\ox K))^*.
   \end{equation}
Finally, if $N_2$ is the kernel of the natural map $V_2\ox \cO\to
E_2$ then
   \begin{equation} \label{obs-ext2}
     {\HH}_{21}^2=H^0(E_1^*\ox N_2\ox K)^*.
   \end{equation}
\end{proposition}

\n {\em Proof. \/} This follows from (\ref{long-exact}) applied to
$(E_1,V_1)=(E,V)$ and $(E_2,V_2)=(E',V')$, together with Serre
duality for the last part. \hfill$\Box$

In order to use this result we will need to be able to estimate the
dimension of ${\HH}_{21}^2$.
\begin{lemma}\label{lem:hopf}
Suppose that $k_2>0$ and $h^0(E_1^*\ox K)\neq 0$. Then the
dimension of $\HH_{21}^2$ is bounded above by
$(k_2-1)(h^0(E_1^*\ox K)-1)$.
\end{lemma}

\n {\em Proof.\/} Use the Hopf lemma which states that, if
$\p:A\ox B \to C$ is a bilinear map between finite-dimensional
complex vector spaces such that, for any $a\in A-\{0\}$,
$\p(a,\cdot)$ is injective and, for any $b\in B-\{0\}$,
$\p(\cdot,b)$ is injective, then the image of $\p$ has dimension
at least $\dim A+\dim B-1$. The result follows from this and
(\ref{obs-ext}). \hfill $\Box$

%%%%%%%%%%%%%%%%%%%%%%%%%%%%%%%%%%%%%%%%%%%%%
\subsection{Infinitesimal deformations}\label{subsec:infinitesimal}
%%%%%%%%%%%%%%%%%%%%%%%%%%%%%%%%%%%%%%%%%%%%%

By standard arguments in deformation theory we have
(see \cite[Th\'eor\`eme 3.12]{He})

\begin{proposition}\label{prop:tangent}
Let $(E,V)$ be an $\a$-stable coherent system.

\begin{enumerate}
\item If $\Ext^2((E,V),(E,V))=0$, then the moduli space of
$\a$-stable coherent systems is smooth in a neighbourhood of the
point defined by $(E,V)$. This condition is satisfied if and only
if the homomorphism $\Ext^1(E,E)\ra \Hom(V,H^1(E))$ is surjective.
\item The Zariski tangent space
to the moduli space at the point defined by $(E,V)$ is
isomorphic to
 $$
 \Ext^1((E,V),(E,V)).
 $$
\end{enumerate}
\end{proposition}

\begin{lemma}\label{lem:BN-number}
Let  $(E,V)$  be an $\a$-stable coherent system of type
$(n,d,k)$. Then
 $$
 \dim \Ext^1((E,V),(E,V))=  \b(n,d,k) + \dim \Ext^2((E,V),(E,V)),
 $$
where $\b(n,d,k)$ is the Brill-Noether number defined in
Definition \ref{def:bnnumber}.
\end{lemma}

\n {\em Proof. \/} By considering  the long exact sequence
(\ref{long-exact}) for $(E',V')=(E,V)$, we see that
\begin{eqnarray}
\dim \Ext^1((E,V),(E,V)) &=&   k\chi (E)-\chi(\End E)- k^2
 + \dim\End(E,V)\nonumber\\
& &\quad\quad \mbox{} +\dim \Ext^2((E,V),(E,V))\nonumber \\
                 &=& k(d+n(1-g))-n^2(1-g)- k^2 +1\nonumber\\
& &\quad\quad\mbox{}+   \dim \Ext^2((E,V),(E,V))\nonumber\\
                 &=& n^2(g-1)+1 -k(k-d +n(g-1))\label{dim-mod}\\
& &\quad\quad\mbox{} + \dim \Ext^2((E,V),(E,V)),\nonumber
\end{eqnarray}
since
 $\End(E,V)\cong\CC$ by Proposition \ref{prop:filtration}(ii). \hfill $\Box$

\begin{corollary}\label{cor:exp}
Every irreducible component $G$ of every moduli space $G_i(n,d,k)$ has
dimension
$$\dim G\ge\b(n,d,k).$$
\end{corollary}

\n{\em Proof. \/} See \cite[Corollaire 3.14]{He}.\hfill$\Box$

The following further corollary of Lemma \ref{lem:BN-number} will be useful.
\begin{corollary}\label{cor:euler}
Let  $C_{21}$ be defined by {\em (\ref{dim-ext})}
and $C_{12}$ by interchanging indices in {\em (\ref{dim-ext})}. Then
 $$
 \b(n,d,k)=\b(n_1,d_1,k_1)+\b(n_2,d_2,k_2)+C_{12}+C_{21}-1.
 $$
\end{corollary}

\n{\em Proof.\/} This follows from (\ref{dim-ext}) and
(\ref{dim-mod}) using the facts that $\chi(E)=\chi(E_1)+\chi(E_2)$
and
 $$
 \chi(\End E)=\chi(\End E_1)+\chi(\End E_2)+\chi(E_1^*\ox
 E_2) +\chi(E_2^*\ox E_1).
 $$
\hfill$\Box$
\begin{remark}\begin{em}\label{rem:duality}
Notice that if $k>n$ then $\b(n,d,k)=\b(k-n,d,k)$. This is easily
seen by writing
 $$
 \b(n,d,k)=n(g-1)(n-k)-k(k-d)+1.
 $$
We will come back  to this symmetry later when studying the dual
span of a coherent system (see section
\ref{subsec:alpha-large.k>n}).
\end{em}\end{remark}
We are now ready to extend to coherent systems the standard fact
about smoothness of Brill-Noether loci. First we extend the definition of
Petri map.
\begin{definition}\label{def:petri}\begin{em}
Let $(E,V)$ be a coherent system. The Petri map of
$(E,V)$ is the map
 $$
 V\ox H^0(E^*\ox K)\lra H^0(\End  E\ox K)
 $$
given by multiplication of sections.\end{em}
\end{definition}

\begin{proposition}\label{prop:Petri}
Let $(E,V)$ be an $\a$-stable coherent system of type
$(n,d,k)$. Then the moduli space $G(\a;n,d,k)$
is smooth of dimension $\b(n,d,k)$ at the point corresponding to
$(E,V)$ if and only if the Petri map is injective.
\end{proposition}

\n {\em Proof. \/} By Proposition \ref{prop:tangent} and Lemma
\ref{lem:BN-number}, the moduli space is smooth of the correct
dimension at $(E,V)$ if and only if $\Ext^2((E,V),(E,V))=0$. The
result is now a special case of (\ref{obs-ext}). \hfill$\Box$

\begin{remark}\begin{em}\label{rem:Petri}
This is a strengthening of the result for Brill-Noether  loci
(Theorem \ref{thm:BN}), and it justifies the idea that the spaces
of coherent systems provide partial desingularisations of the
Brill-Noether loci (see sections \ref{subsec:alpha-small} and
\ref{sec:applications}). In view of Proposition \ref{prop:Petri} and
Corollary \ref{cor:exp}, we often refer to $\b(n,d,k)$ as the {\em expected
dimension} of $G_i(n,d,k)$.
\end{em}\end{remark}
There is a special case in which it is easy to check the injectivity of the
Petri map.

\begin{proposition}\label{prop:smooth}
Let $(E,V)$ be an $\a$-stable coherent system such that  $k\leq
n$. If $V\ox \cO\to E$ is injective then the moduli space is smooth
of dimension $\b(n,d,k)$ at the point corresponding to $(E,V)$.
This happens in particular when $k=1$.
\end{proposition}
\n {\em Proof. \/} We have an exact sequence
 $$
  0\lra V\ox\cO\lra E\lra F\lra 0,
  $$
where $F$ is a coherent sheaf. Tensoring with $E^*\ox K$ gives
 $$
 0\lra V\ox E^*\ox K\lra \End E\ox K\lra F\ox E^*\ox K\lra0.
 $$
Taking sections, we see that the Petri map is
injective.\hfill$\Box$

When the hypotheses of Proposition \ref{prop:smooth} fail, we will need to
estimate
the dimension of the kernel of the Petri map. In fact Lemma
\ref{lem:hopf} gives us such an estimate.

%%%%%%%%%%%%%%%%%%%%%%%%%%%%
\section{Range for $\a$}\label{sec:range}
%%%%%%%%%%%%%%%%%%%%%%%%%%%%

\begin{lemma}\label{lem:range}
If $k<n$ then the moduli space of $\a$-semistable coherent systems
of type $(n,d,k)$ is empty for  $\a> \frac{d}{n-k}$. In
particular, we must have  $d\ge 0$ in order for $\a$-semistable
coherent systems to exist. Also we must have $d>0$ in order  for
$\a$-stable coherent systems to exist.
\end{lemma}
\n {\em Proof. \/} Suppose that  $(E,V)$ is an $\a$-semistable
coherent system of type $(n,d,k)$. By applying the
$\a$-semistability condition to the subsystem $(E',V)$, where
$E'=\Im(V\ox\cO\ra E)$, one obtains the upper bound $\a\leq
\frac{d}{n-k}$, which in particular implies that $d\geq 0$ in
order to have non-empty moduli spaces. The final assertion is
similar. \hfill $\Box$

{}From this lemma and the considerations of section \ref{subsec:coh},
we deduce at once the following proposition.
\begin{proposition}\label{prop:last}
Let $k<n$ and let $\a_L$ be  the biggest critical value smaller
than $\frac{d}{n-k}$. The $\a$-range is then divided in a finite
set of intervals determined by
 $$
 0=\a_0<\a_1<\dots <\a_L<\frac{d}{n-k}\ .
 $$
Moreover, if $\a_i$ and $\a_{i+1}$ are two consecutive critical
values, the moduli spaces for two different values of $\a$ in the
interval $(\a_i,\a_{i+1})$ coincide, and if $\a>
\frac{d}{n-k}$ the moduli spaces are empty.
\end{proposition}

\begin{lemma}\label{lem:degree}
Let $k\geq n$. We must have  $d\ge 0$ in order for $\a$-semistable
coherent systems of type $(n,d,k)$ to exist. Also we must have
$d>0$ in order for $\a$-stable coherent systems to exist except in
the case $(n,d,k)=(1,0,1)$.
\end{lemma}
\n {\em Proof. \/} The first assertion is clear if the  map
$\cO\ox V\ra E$ is generically surjective, otherwise one has to
apply the $\a$-semistability condition to the subsystem $(E',V)$,
where $E'=\Im(\cO\ox V\ra E)$.

For the second assertion, suppose $d=0$ and apply the $\a$-stability
condition to
$(E',V)$, where $E'=\Im(\cO\ox V\ra E)$, to get that $\cO\ox V\ra
E$ is generically surjective. Therefore $E\cong \cO^n$ and
$\a$-stability forces $k=n=1$. \hfill $\Box$

Although  in the case $k\geq n$ the  stability condition does not
provide us with a bound for $\a$, we will show that in fact after a
certain finite value of $\a$ the moduli spaces do not change. We
show first that for $\a$ big enough the vector bundle $E$ for an
$\a$-semistable coherent system is generically generated by the
sections in $V$. More precisely
\begin{proposition}\label{prop:BGN-last}
Suppose $k\geq n$. Then there exists $\a_{gg}>0$ such that for
$\a\geq\a_{gg}$ if $(E,V)$ is $\a$-semistable then the map
$\p:V\ox \cO\ra E$ is generically surjective, i.e. we have an
exact sequence
  $$
  0\lra N\lra\cO^{\oplus k}\stackrel{\p}{\lra}
  E\lra T\lra 0\
  $$
where
\begin{enumerate}
 \item $T$ is a torsion sheaf (possibly $0$),
 \item $\rk N=k-n$,
 \item $H^0(N)= 0$.
\end{enumerate}
In fact,
  $$
  \a_{gg}\leq \frac{d(n-1)}{k}.
  $$
\end{proposition}
\n {\em Proof. \/} Let $N=\Ker \p$ and $I=\Im \p$, and suppose
$\rk I=n-l<n$. One has  the exact sequence
 $$
 0\lra N\lra\cO^{\oplus k}\lra E\lra E/I\lra 0.
 $$
Consider the subsystem $(I,V)$. One has $d_I:=\deg I\geq 0$ since
$I$ is generated by global sections. Now $\a$-semistability
implies that $\mu_\a(I,V)\leq\mu_\a(E,V)$, which means that $$
\frac{d_I}{n-l}+\a\frac{k}{n-l} \leq\frac{d}{n} +\a\frac{k}{n}, $$
and hence
 $$
  \a\leq \frac{d(n-l)}{kl}\leq\frac{d(n-1)}{k}.
 $$
We conclude that if $\a>\frac{d(n-1)}{k}$ then $\rk I=n$ so that
$E/I=T$ is pure torsion. Finally $H^0(N)=0$ since $H^0(\p)$ is
injective by definition of coherent system. \hfill $\Box$

Our next  object  is to show that the $\a$-stability condition is
independent of $\a$ for $\a>d(n-1)$. More precisely
\begin{proposition}\label{stabilization}
\begin{enumerate}
 \item If there exists a subsystem $(E',V')$ of $(E,V)$ with
$\displaystyle{\frac{k'}{n'}>\frac{k}{n}}$, then $(E,V)$ is not
$\a$-semistable for $\a>d(n-1)$.
 \item If there exists a subsystem $(E',V')$ of $(E,V)$ with
$\displaystyle{\frac{k'}{n'}=\frac{k}{n}\hbox{ and
}\frac{d'}{n'}\ge\frac{d}{n}}$, then $(E,V)$ is not $\a$-stable
for any $\a$.
 \item If neither {\em (i)} nor {\em (ii)} holds, and $E$ is generically
generated by its sections, then $(E,V)$ is $\a$-stable for
$\a>d(n-1)$.
 \end{enumerate}
\end{proposition}

\n{\it Proof}. (i) Suppose $(E,V)$ is $\a$-semistable. Replacing
$E'$ by a subbundle if necessary, we can suppose that $E'$ is
generically generated by its sections and hence $d'\ge0$. Then we
have
 $$
 \a\frac{k'}{n'}\le \frac{d}{n}+\a\frac{k}{n},
 $$
i.e.
 $$
 \a\le\frac{n'd}{nk'-n'k}\le d(n-1).
 $$

(ii) is obvious.

(iii) If neither (i) nor (ii) holds and $(E',V')$ contradicts the
$\a$-stability of $(E,V)$, then we must have
$\displaystyle{\frac{k'}{n'}<\frac{k}{n}}$. If $E$ is generically
generated by its sections, then so is $E/E'$; hence
$\deg(E/E')\ge0$ and $d'=\deg E'\le d$. Thus we have
 $$
 \frac{d}{n}+\a\frac{k}{n}\le \frac{d'}{n'}+\a\frac{k'}{n'}\le
 \frac{d}{n'}+\a\frac{k'}{n'},
 $$
i.e.
 $$
 \a\le\frac{d(n-n')}{n'k-nk'}\le d(n-1).
 $$
\hfill $\Box$

We have thus proved the following.

\begin{proposition}\label{prop:flips}
Let $k\geq n$. Then there is a critical value, denoted by $\a_L$,
after which the moduli spaces stabilise, i.e. $G(\a)=G_L$ if
$\a>\a_L$. The $\a$-range is thus divided into a finite set of
intervals bounded by critical values
 $$
 0=\a_0<\a_1<\dots<\a_L<\infty
 $$
and such that
\begin{enumerate}
\item if $\a_i$ and $\a_{i+1}$ are two consecutive
critical values, the moduli spaces for any two different values
of $\a$ in the interval $(\a_i,\a_{i+1})$ coincide,
\item
for any two different values of $\a$ in the range $(\a_L,\infty)$,
the moduli spaces coincide.
\end{enumerate}
\end{proposition}

%%%%%%%%%%%%%%%%%%%%%%%%%%%%%%%%%%%%%%%%%%%%%%%%%%%%%%%%%%%%%%%%%%%%
\section{Moduli for $\a$ large}\label{sec:alpha-large}
%%%%%%%%%%%%%%%%%%%%%%%%%%%%%%%%%%%%%%%%%%%%%%%%%%%%%%%%%%%%%%%%%%%

%%%%%%%%%%%%%%%%%%%%%%%%%%%%%%%%%%%%%%%%%%%%%%%%%%%%%%%%%%%%%
\subsection{The moduli space $G_L$ for  $k<n$}\label{subsec:alpha-large.k<n}
%%%%%%%%%%%%%%%%%%%%%%%%%%%%%%%%%%%%%%%%%%%%%%%%%%%%%%%%%%%%%
Recall that, when $k<n$,  $G_L$ denotes the moduli space of coherent systems
for $\a$ large, i.e. $\a_L<\a<\frac{d}{n-k}$.
The description of $G_L$ in this case has been carried out
in \cite{BG2}, where we refer for details. We summarise here the
main results.

\begin{definition}\label{def:BGN}\begin{em}
A {\em BGN extension} \textup{(}\cite{BGN}\textup{)} is an extension of
vector bundles
 $$
 0\lra \cO^{\oplus k}\lra E\lra  F\lra 0
 $$
which satisfies the following conditions:
\begin{itemize}
\item $\rk E=n>k$,
\item $\deg E=d>0$,
\item $H^0(F^*)=0$,
\item if $\vec{e}=(e_1,\dots,e_k)\in H^1(F^\ast\ox
\cO^{\oplus k})=H^1(F^\ast)^{\oplus k}$\ denotes the
class of the extension, then  $e_1,\dots,e_k$\ are linearly independent as
vectors in $H^1(F^\ast)$.
\end{itemize}\end{em}
\end{definition}
The BGN extensions which differ only by an automorphism of
$\cO^{\oplus k}$, i.e.~by the action of an element in ${\GL}(k)$,
comprise a  BGN extension class of type $(n,d,k)$.

\begin{proposition} \label{prop:BGN}
Suppose that $0<k<n$\ and $d>0$.
Let $\a_L<\a <\frac{d}{n-k}$. Let $(E,V)$ be an $\a$-semistable
coherent system of type $(n,d,k)$. Then $(E,V)$ defines a BGN
extension class represented by an extension
 \begin{equation}\label{bgn}
 0\lra \cO^{\oplus k}\lra E\lra F\lra 0,
  \end{equation}
with $F$ semistable. In the converse direction,  any BGN extension
of type $(n,d,k)$ in which the quotient $F$\ is stable gives rise
to an $\a$-stable coherent system of the same type.
\end{proposition}
\begin{remark}\begin{em} In the last part of Proposition \ref{prop:BGN},
it is essential to have $F$ stable. If $F$ is only semistable,
the coherent system can fail to be $\a$-semistable.
\end{em}\end{remark}

\begin{theorem}\label{G_L(k<n)}
Let $0<k<n$\ and $d>0$. If $g\ge2$, the moduli space $G_L(n,d,k)$
of $\a$-stable coherent systems of type $(n,d,k)$ is birationally
equivalent to a fibration over the moduli space $M(n-k,d)$ of
stable bundles of rank $n-k$ and degree $d$ with fibre the
Grassmannian ${\Gr}(k,d+(n-k)(g-1))$. In particular $G_L$ is
non-empty if and only if $k\leq d+(n-k)(g-1)$, and it is then
always irreducible and smooth of dimension $\b(n,d,k)$. If
$(n-k,d)=1$ then the birational equivalence is an isomorphism.
\end{theorem}

\n {\em Proof. \/} This follows directly from Proposition \ref{prop:BGN}
and the remark immediately preceding it. \hfill$\Box$

\begin{remark}\label{rmk:g=0,1}\begin{em}
Proposition \ref{prop:BGN} remains true when $g=0$ or $1$, but Theorem
\ref{G_L(k<n)} can fail because $M(n-k,d)$ may be empty. In fact, if $g=0$,
$M(n-k,d)=\emptyset$ unless $n-k=1$. Furthermore, if $n-k$ does not divide
$d$, then $\tilde{M}(n-k,d)=\emptyset$, and it follows from Proposition
\ref{prop:BGN} that $G_L(n,d,k)=\emptyset$. If $d=(n-k)a$ with $a\in\ZZ$,
then $\tilde{M}(n-k,d)$ consists of a single point corresponding to the
bundle
$$F={\mathcal O}(a)\oplus\ldots\oplus{\mathcal O}(a).$$
It is not clear from the results of \cite{BG2} whether there exist
$\alpha$-stable coherent systems as in (\ref{bgn}). Thus we conclude,
for $g=0$,
\begin{itemize}
\item $G_L(n,d,n-1)\ne\emptyset$ if and only if $d\ge n$, and it is then
isomorphic to the Grassmannian ${\Gr}(n-1,d-1)$,
\item $G_L(n,d,k)=\emptyset$ if $k\le n-2$ and $d$ is not divisible by
$n-k$,
\item if $k\le n-2$ and $d=(n-k)a$ with $a\in\ZZ$, then
$G_L(n,d,k)=\emptyset$ if $d<n$; if $d\ge n$, a more detailed analysis is
required.
\end{itemize}

Turning now to $g=1$, we know that $\tilde{M}(n-k,d)$ is always non-empty
and that $M(n-k,d)\ne\emptyset$ if and only if $(n-k,d)=1$; moreover in this
case $M(n-k,d)$ is isomorphic to the curve $X$. (All this follows
essentially
from \cite{At}.) We conclude, for $g=1$,
\begin{itemize}
\item if $(n-k,d)=1$, then $G_L(n,d,k)\ne\emptyset$ if and only if $d\ge k$,
and it is then isomorphic to a fibration over $X$ with fibre
${\Gr}(k,d)$,
\item if $(n-k,d)\ne1$, a more detailed analysis is required.
\end{itemize}
\end{em}\end{remark}

%%%%%%%%%%%%%%%%%%%%%%%%%%%%%%%%%%%%%%%%%%%%
\subsection{Quot schemes}\label{subsec:quot}
%%%%%%%%%%%%%%%%%%%%%%%%%%%%%%%%%%%%%%%%%%%%
When $k\ge n$, we can follow \cite{BeDW} and relate $G_L$ to a
Grothen\-dieck
Quot scheme. In fact, by Proposition \ref{prop:BGN-last}, any element
of $\tilde{G}_L$ can be represented in the form
\begin{equation}\label{quot}
0\lra N\lra\VV\ox\cO\stackrel{\p}{\lra}E,
\end{equation}
where $\p$ is generically surjective. Dualising (\ref{quot}), we obtain
\begin{equation}\label{quot2}
0\lra E^*\lra\VV^*\ox\cO\lra F\lra0,
\end{equation}
where $F$ is a coherent sheaf but is not torsion-free (unless $\p$
is surjective). Conversely, given (\ref{quot2}), one can recover
(\ref{quot}) (in fact $N\cong F^*$). It follows that there is a
bijective correspondence between triples $(E,\VV,\p)$ and points
of $Q=\hbox{Quot}_{k-n,d}(\cO^{\oplus k})$, the Quot scheme of
quotients of $\cO^{\oplus k}$ of rank $k-n$ and degree $d$. In
order to obtain $\tilde{G}_L$, we therefore need to construct a
GIT quotient of $Q$ by the natural action of ${\GL}(k)$ with
respect to a stability condition corresponding to the
$\a$-stability of coherent systems for large $\a$. This situation
requires detailed analysis, but even if we complete the
construction, it may still be difficult to obtain information
about $G_L$, since even basic information about $Q$ is often
lacking, e.g., when it is non-empty, irreducible etc. However,
potentially this would be a useful source of information about
$G_L$.

In sections \ref{subsec:alpha-large.k=n} and
\ref{subsec:alpha-large.k>n}, we shall use the sequences
(\ref{quot}) and (\ref{quot2}) to obtain information about $G_L$
in the cases $k=n$ and $k>n$.

%%%%%%%%%%%%%%%%%%%%%%%%%%%%%%%%%%%%%%%%%%%%%
\subsection{The moduli space $G_L$ for  $k=n$}\label{subsec:alpha-large.k=n}
%%%%%%%%%%%%%%%%%%%%%%%%%%%%%%%%%%%%%%%%%%%%%
We are now able to prove Theorem \ref{thm:main1} in a stronger form which
covers
$\tilde{G}_L$ as well as $G_L$.

\begin{theorem}\label{G_L(k=n)}
Let $k=n\ge2$.  If $\tilde{G}_L$ is non-empty then $d=0$ or $d\geq n$.
For $d>n$,
  $\tilde{G}_L$ is irreducible and $G_L$ is smooth  of the expected
dimension $dn-n^2+1$.
  For $d=n$, $G_L$ is empty and $\tilde{G}_L$ is irreducible
of dimension $n$ (not of the expected dimension). For $d=0$,
$G_L$ is empty and $\tilde{G}_L$ consists of a single point.
\end{theorem}

\n {\em Proof. \/} Let $(E,V)$ be an $\a$-semistable coherent
system for any $\a$, represented by $\p:\cO^{\oplus n}\ra E$. If
$\p$ is an isomorphism, then $d=0$ and $(E,V)\cong(\cO^{\oplus
n},\CC^n)$, which is clearly $\a$-semistable but not $\a$-stable.
Otherwise there exists $\cO\subset\cO^{\oplus n}$ which defines a
section of $E$ with a zero. It follows that this section is
contained in a subbundle of $E$ of rank $1$ with degree $>0$. This
subbundle together with the section defines a subsystem which
contradicts $\a$-semistability if $d<n$ and $\a$-stability if
$d=n$.

Now suppose $d\ge n$.
  Let $(E,V)$ be any $\a$-semistable coherent system for $\a$
  large. By Proposition \ref{prop:BGN-last} with $k=n$, we have an extension
 $$
  0 \to {\cO}^{\oplus n} \to E \to T \to 0,
  $$
where $T$ is torsion.
  The generic torsion sheaf is of the form $T={\cO}_D$
  for a divisor $D$ consisting of $d$ distinct points.
For such $T$, $E$ is given
  by an extension class $\xi \in \Ext^1(T,{\cO}^{\oplus n}) \cong
  \Hom (T, {\CC}^n)$,
  which is equivalent to a collection of $d$ vectors $\xi_i\in
  {\CC}^n$, one for each point $P_i$ in the support of $D$.

  We claim that all the coherent systems defined by extensions in
 $$
  U=\{ ({\cO}_D, \xi)| \text{any subset of
  $n$ vectors of $\xi_1,\ldots, \xi_d$ is linearly independent}\}
 $$
  are $\a$-stable (for $d>n$) or $\a$-semistable
  (for $d=n$).

Suppose for the moment that the claim holds. Let
 $$
  U^{ss}=\{ (T, \xi)\;|\; \xi\in \Ext^1(T,\cO^{\oplus n}) \text{ and
 determines an $\a$-semistable  coherent system}\}.
 $$
  Then
  $U$ is dense and open in $U^{ss}$. Also $U^{ss}$
  dominates the moduli space of $\a$-semistable coherent
  systems. Thus $\tilde{G}_L$ is irreducible. The fact that $G_L$
  is smooth of the expected dimension follows at once from
  Proposition \ref{prop:smooth}. We can also compute directly
  the dimension of the space of coherent systems determined by $U$.
  The space $\Ext^1(T,{\cO}^{\oplus n})$ has dimension $dn$, and we have
  to quotient out by the automorphisms $\Aut_{\cO} T=
  \GL(1)^d$, for $T={\cO}_D$, and by
  $\Aut_{\cO} {\cO}^{\oplus n}=\GL(n)$. For $d>n$,
  the centraliser of the action of the product on $(\cO_D,\xi) \in U$ is
  ${\CC}^*$. So the dimension of the space of coherent systems
  determined by $U$ is $d+dn-d-n^2+1$, which is in agreement with
  our previous answer.

  For $d=n$, $\GL(n)$ acts freely on
  any collections of $n$ linearly independent
  vectors in ${\CC}^n$. So the dimension of the space of coherent
  systems determined by $U$ is $d+dn-n^2=n$. It is possible
  that different elements of $U$ give rise to S-equivalent systems,
  thus reducing $\dim \tilde{G}_L$. However, if $g\ge1$, the coherent
systems
  $$
  \left(\cO(P_1),H^0(\cO(P_1))\right)\oplus\ldots
  \oplus\left(\cO(P_n),H^0(\cO(P_n))\right),
  $$
where $P_1,\ldots,P_n\in X$, are clearly $\a$-semistable and no two of them
are S-equivalent, so $\dim\tilde{G}_L\ge n$, which completes the
computation. If $g=0$, there is a unique line bundle $\cO(1)$ of degree $1$,
and $h^0(\cO(1))=2$; in this case the coherent systems
$$(\cO(1),V_1)\oplus\ldots\oplus(\cO(1),V_n),$$
where $V_1,\ldots,V_n$ are subspaces of dimension $1$ of $H^0(\cO(1))$,
are $\a$-semistable and form
a family of dimension $n$. This gives the same conclusion.

  Now we prove the claim, i.e.\ every $(E,V)$ in
  the image of $U$ is $\a$-stable for
  $\a$ large. Let $(E_1,V_1)$ be a coherent subsystem of $(E,V)$. As
  $E_1 \subset E$ we must have $k_1 \leq n_1$. If $k_1<n_1$ then $(E_1,V_1)$
  cannot violate $\a$-stability. If $k_1=n_1$ then we have a diagram
 $$
 \begin{array}{ccccc}
 {\cO}^{\oplus n_1} & \to & E_1 &\to & T_1 \\
   \downarrow & & \downarrow & & \downarrow \\
 {\cO}^{\oplus n} & \to & E &\to & T
 \end{array}
 $$
  where $d_1=\deg T_1$. Then the image of $\xi \in \Ext^1(T, {\cO}^{\oplus
n})$
  in $\Ext^1(T_1, {\cO}^{\oplus n})$ lies in the subspace
  $\Ext^1(T_1, {\cO}^{\oplus n_1})$. This is equivalent to
  $\xi_i \in {\CC}^{n_1}$ for any $P_i$ in the support of $T_1$.
  $E$ is $\a$-semistable if $d_1/n_1 \leq d/n$ for all possible
  choices of diagrams as above.

  Now any subcollection of $n$ vectors of the $\xi_i$ is linearly
  independent, so for $d_1\geq n$ we have $n_1=n$ and $E_1=E$.
For $d_1<n$ we have $n_1\ge d_1$ and $d_1/n_1\le 1$.
Hence, for $d>n$ the coherent systems are $\a$-stable, while
  for $d=n$ they are $\a$-semistable.\hfill$\Box$

\begin{remark}\begin{em}\label{rem:d=n}
For $d\le n$, the proof shows that $(E,V)$ cannot be $\a$-stable
for any $\a$; moreover, if $0<d<n$, $(E,V)$ cannot be
$\a$-semistable. For $d=0$, any $\a$-semistable coherent system is
isomorphic to $(\cO^{\oplus n},\CC^n)$. For $d=n$, one can show
that $\tilde{G}(\a)$ is independent of $\a$ and that
$\tilde{G}(\a)\cong S^nX$ \textup{(}see \cite[Theorem 8.3]{BGN}
for the case $\a=0$\textup{)}.
\end{em}\end{remark}

%%%%%%%%%%%%%%%%%%%%%%%%%%%%%%%%%%%%%%%%%%%%%%%%%%%%%%%%%%%%%%%
\subsection{The moduli space $G_L$ for  $k> n$. The dual span construction}
\label{subsec:alpha-large.k>n}
%%%%%%%%%%%%%%%%%%%%%%%%%%%%%%%%%%%%%%%%%%%%%%%%%%%%%%%%%%%%%%%
We can represent a coherent system by a sequence (\ref{quot}),
where we now suppose that $k>n$ and that $\p$ is surjective; so we have
 $$
 0\lra N\lra\VV\ox\cO\stackrel{\p}{\lra}E\lra0
 $$
and
 $$
 0\lra E^*\lra\VV^*\ox\cO\stackrel{\psi}{\lra} N^*\lra0.
 $$
In the case where
$V=H^0(E)$ and $E$ is generated by its sections, this construction
has been used by a number of authors (see for example
\cite{Bu1,Bu2,EL,M1,PR}),
the main question being to determine conditions under
which the stability of $E$ implies that of $N$. Recently Butler
noted that the construction belongs more naturally to the theory
of coherent systems and began to investigate it using
$\a$-stability. However he restricted attention to small $\a$. Our
purpose in this section is to show that the construction works
better if we consider large $\a$.

It is convenient here to make partial use of the wider notion of coherent
system,
introduced in \cite{KN} and mentioned in section \ref{subsec:coh},
by dropping the assumption that $H^0(\p)$ is injective. This
makes no essential difference as $(E,V)$ cannot be
$\a$-semistable unless $H^0(\p)$ is injective.
It does however mean that $(N^*,\VV^*,\psi)$ always determines a coherent
system, which we may call the {\it dual span} of
$(E,\VV,\p)$ and denote by $D(E,\VV,\p)$ (or $D(E,V)$ in the case where
$H^0(\p)=0$).

\begin{definition}\label{def:strongly-unstable}\begin{em}
$(E,\VV,\p)$ is {\em strongly unstable} if there exists a proper
coherent subsystem $(E',\VV',\p')$ such that
 $$
 {\frac{k'}{n'}>\frac{k}{n}}\, .
 $$\end{em}
\end{definition}

\begin{proposition}\label{prop:dual-strongly}
Suppose that $E$ is generated by $\VV$. Then
$(E,\VV,\p)$ is strongly unstable if and only if
$D(E,\VV,\p)$ is strongly unstable.
\end{proposition}

\n{\it Proof}. Suppose that $(E,\VV,\p)$ is strongly unstable and
that $(E',\VV',\p')$ is a subsystem as in the definition above.
Replacing $E'$ by the (sheaf) image of $\VV'$ in $E$ if necessary,
we have a short exact sequence
 $$
 0\to N'\to\VV'\ox\cO\to E'\to0.
 $$
$D(E',\VV',\p')$ is then a quotient system of $D(E,\VV,\p)$. The
corresponding subsystem has rank $(k-n)-(k'-n')$ and dimension
$k-k'$ and
 $$
 (k-n)(k-k')-((k-n)-(k'-n'))k=nk'-n'k>0.
 $$
So $D(E,\VV,\p)$ is strongly unstable. The converse is similar,
which completes the proof. \hfill $\Box$

By Proposition \ref{stabilization}, any strongly unstable coherent
system fails to be $\a$-semistable for $\a>d(n-1)$. The converse
may fail because we have to take account of subsystems with
$\displaystyle{\frac{k'}{n'}=\frac{k}{n}}$. However, if $(n,k)=1$,
there are no such subsystems and we have
\begin{corollary}\label{cor:dual-span}
Suppose that $E$ is generated by $V$. If $(n,k)=1$, then
$(E,V)$ is $\a$-stable for large $\a$ if and only if $D(E,V)$
is $\a$-stable for large $\a$.
\end{corollary}

By Proposition \ref{stabilization} it is sufficient to take
$\a>\max\{d(n-1),d(k-n-1)\}$.

\begin{theorem}\label{thm:dual-span}
Suppose that $X$ is a Petri curve and that $k=n+1$. Then $G_L$ is non-empty
if and only if $\b=g-(n+1)(n-d+g)\ge0$. Moreover $G_L$ has dimension $\b$
and it is
irreducible whenever $\b>0$.

\end{theorem}

{\em Proof.\ } If $(E,V)\in G_L$, then by Proposition
\ref{prop:BGN-last}, $E$ is generically generated by $V$. If we
suppose further that $E$ is generated by $V$, then $D(E,V)\in
G(1,d,n+1)$. Since $X$ is Petri,
 $G(1,d,n+1)$ is non-empty if
and only if the Brill-Noether number
 $$
 \b=\b(1,d,n+1)=g-(n+1)(n-d+g)\ge0.
 $$
Moreover, if this holds, $G(1,d,n+1)$ has dimension $\b$, and it
is irreducible whenever $\b>0$. Note also that the dimension of
the subvariety consisting of systems $(L,W)$ for which $L$ is not
generated by $W$ has dimension at most
 $$
 g-(n+1)(n-(d-1)+g)+1<\b.
 $$
So $G(1,d,n+1)$ has a dense open subset in which $L$ is generated
by $W$. The Brill-Noether number $\b(n,d,n+1)=\b$ by Remark
\ref{rem:duality} so the systems $(E,V)$ which are $\a$-stable for
large $\a$ and for which $V$ generates $E$ are parametrised by a
variety of the expected dimension. If $E$ is only generically
generated by $V$ and $E'$ is the subsheaf generated by $V$, we can
put $\deg E'=d-t$ with $t>0$. Then, by the argument above, the
variety parametrising the systems $(E',V)$ has the expected
dimension, which is $\b-(n+1)t$. On the other hand, the variety
parametrising the extensions
 $$
 0\to E'\to E\to T\to 0,
 $$
where $T$ is a torsion sheaf of length $t$, has dimension $nt$
(after factoring out by the action of $\Aut T$). So the variety
parametrising all the corresponding systems $(E,V)$ has dimension
$<\b$. Since every component of $G_L$ has dimension $\ge\b$, this
completes the proof. \hfill$\Box$

%%%%%%%%%%%%%%%%%%%%%%%%%%%%%%%%%%%%%
\section{Crossing critical values}\label{sec:critical}
%%%%%%%%%%%%%%%%%%%%%%%%%%%%%%%%%%%%%

In this section we analyse the differences between consecutive
moduli spaces in the family $\{G_0,G_1,\dots,G_L\}$. Recall that
$G_i$ denotes the moduli space of $\a$-stable coherent systems,
where $\a$ is (anywhere) in the interval bounded by the critical
values $\a_i$ and $\a_{i+1}$. The differences between $G_{i-1}$
and $G_i$ are thus due to the differences between the
$\a$-stability conditions for $\a<\a_i$ and $\a>\a_i$.

%%%%%%%%%%%%%%%%%%
\subsection{The basic mechanism}\label{subsec:mechanism}
%%%%%%%%%%%%%%%%%%%

The following lemma describes the basic mechanism responsible for a
change in the stability property of a coherent system.

\begin{lemma}
Let $(E,V)$ be a coherent system of type $(n,d,k)$ and let
$(E',V')$ be  a subsystem of type $(n',d',k')$. Then
$\mu_{\a}(E',V')-\mu_{\a}(E,V)$ is a linear function of $\a$ which
is

\begin{itemize}
\item monotonically increasing if  $\frac{k'}{n'}- \frac{k}{n}>0$,
\item monotonically decreasing if  $\frac{k'}{n'}- \frac{k}{n}<0$,
\item constant if $\frac{k'}{n'}- \frac{k}{n}=0$.
\end{itemize}

In particular, if $\a_i$ is a critical value and
$\mu_{\a_i}(E',V')=\mu_{\a_i}(E,V)$, then
\begin{itemize}
 \item $(\mu_{\a}(E',V')-\mu_{\a}(E,V))(\a-\a_i)>0$, for all
 $\a\ne\a_i$, if $\frac{k'}{n'}- \frac{k}{n}>0$,
 \item $(\mu_{\a}(E',V')-\mu_{\a}(E,V))(\a-\a_i)<0$, for all
 $\a\ne\a_i$, if $\frac{k'}{n'}- \frac{k}{n}<0$,
 \item $\mu_{\a}(E',V')-\mu_{\a}(E,V)=0$ for all
 $\a$ if $\frac{k'}{n'}- \frac{k}{n}=0$.
\end{itemize}
\end{lemma}
\n {\em Proof. \/} This follows easily from
 $$
 \mu_{\a}(E',V')-\mu_{\a}(E,V)=\frac{d'}{n'}-\frac{d}{n}+
 \a \left(\frac{k'}{n'}- \frac{k}{n}\right).
 $$
\hfill $\Box$

In particular, we have
\begin{lemma}  Let $( E,V)$\ be a coherent system of type
$\ndk$. Suppose that it is $\a$-stable for  $\a>\a_i$, but is
strictly $\a$-semistable for $\a=\a_i$.  Then $( E,V)$\ is
unstable for all $\a<\a_i$.
\end{lemma}

\n {\em Proof. \/} Any such coherent system must have a subsystem,
say $(E',V')$, for which $\mu_{\a_i}(E',V')=\mu_{\a_i}(E,V)$ but
such that $\mu_{\a}(E',V')<\mu_{\a}(E,V)$ if $\a>\a_i$. It follows
from the previous lemma that $\mu_{\a}(E',V')>\mu_{\a}(E,V)$ for
all $\a<\a_i$, i.e. the subsystem $(E',V')$ is destabilising for
all $\a<\a_i$. \qed

\begin{remark}\begin{em} Thus, if we study the effect on
$G_L\ndk$\ of monotonically reducing $\a$, we see that ``once a
coherent system is removed it can never return''. In contrast to
this, it can happen that ``a coherent system once added may have
to be later removed'' \cite{BG2}.
\end{em}\end{remark}

\begin{definition}\begin{em}
We define $G_i^+\subseteq G_i=G_i(n,d,k)$ to be the set of all
$(E,V)$ in $G_i$ which are not $\a$-stable if $\a<\a_i$.
Similarly, we define $G_i^-\subseteq G_{i-1}$ to be the set of all
$(E,V)$ in $G_{i-1}$ which are not $\a$-stable if $\a>\a_i$.\end{em}
\end{definition}

We can identify the sets $G_i-G_i^+ = G_{i-1}-G_i^- $ and hence
(set theoretically) we get

\begin{itemize}
\item $G_{i+1}=G_i-G_{i+1}^- +G_{i+1}^+$,
\item $G_{i-1}=G_i-G_i^+ +G_i^-$.
\end{itemize}

\noindent In fact, we can be more precise.  The subset $G_i^+$ consists
of the points in $G_i$ corresponding to coherent systems which are not
$\a_i$-stable; they therefore form a closed subscheme of $G_i$. Similarly
$G_i^-$ is a closed subscheme of $G_{i-1}$. Hence $G_i-G_i^+$ and
$G_{i-1}-G_i^-$ have natural scheme structures, and as such are isomorphic.
%%%%%%%%%%%%%%%%%%%
\subsection{Destabilising patterns}\label{subsec:destab}
%%%%%%%%%%%%%%%%%%%

The following lemma allows us to describe the sets $G_i^+$ and
$G_i^-$, and also to estimate their codimensions in the moduli
spaces $G_i$. It is important to note that, unlike the
Jordan-H\"older filtrations for semistable objects, the
descriptions we obtain are always as extensions, i.e. 1-step
filtrations. This simplification results from a careful
exploitation of the stability parameter. For convenience, we
denote values of $\a$ in the intervals on either side of $\a_i$ by
$\a_i^-$ and $\a_i^+$ respectively.

\begin{lemma}\label{lem:vicente} Let $\a_i$ be a critical value
of $\a$ with $1\le i\le L$. Let $(E,V)$ be a coherent system of
type $\ndk$.

\begin{enumerate}
\item Suppose that $(E,V)$ is
$\a_i^+$-stable but $\a^-_i$-unstable.  Then $(E,V)$ appears as
the middle term in an extension
 \begin{equation}\label{destab}
 0\to (E_1,V_1)\to (E,V)\to (E_2,V_2) \to 0
 \end{equation}
in which
 \begin{enumerate}
 \item  $(E_1,V_1)$ and $(E_2,V_2)$ are both $\a^+_i$-stable, with
 $\mu_{\a_i^+}(E_1,V_1)<\mu_{\a_i^+}(E_2,V_2)$,
 \item $(E_1,V_1)$ and $(E_2,V_2)$ are both $\a_i$-semistable,
 with $\mu_{\a_i}(E_1,V_1)=\mu_{\a_i}(E_2,V_2)$,
 \item $\frac{k_1}{n_1}$ is a maximum among all proper subsystems
 $(E_1,V_1)\subset (E,V)$  which satisfy {\rm (b)},
 \item $n_1$ is a minimum among all subsystems which satisfy {\rm (c)}.
 \end{enumerate}

\item Similarly, if $(E,V)$ is $\a_i^-$-stable but
$\a^+_i$-unstable, then $(E,V)$ appears as the middle term in an
extension $\mathrm{(\ref{destab})}$ in which
 \begin{enumerate}
 \item $(E_1,V_1)$ and $(E_2,V_2)$ are both $\a^-_i$-stable, with
 $\mu_{\a_i^-}(E_1,V_1)<\mu_{\a_i^-}(E_2,V_2)$,
 \item $(E_1,V_1)$ and $(E_2,V_2)$ are both $\a_i$-semistable,
 with $\mu_{\a_i}(E_1,V_1)=\mu_{\a_i}(E_2,V_2)$,
 \item $\frac{k_1}{n_1}$ is a minimum among all proper subsystems
 $(E_1,V_1)\subset (E,V)$  which satisfy {\rm (b)},
 \item $n_1$ is a minimum among all subsystems which satisfy {\rm (c)}.
 \end{enumerate}
\end{enumerate}
\end{lemma}

\n {\em Proof.\/} Since its stability property changes at $\a_i$,
the coherent system $(E,V)$ must be strictly $\a_i$-semistable,
i.e. it must have a proper subsystem $(E',V')$ with
$\mu_{\a_i}(E',V')=\mu_{\a_i}(E,V)$. Consider the (non-empty) set
 $$
 {\cF}_1=\{(E_1,V_1)\subsetneq (E,V)\ |\
 \mu_{\a_i}(E_1,V_1)= \mu_{\a_i}(E,V)\ \}.
 $$
Any such subsystem $(E_1,V_1)$ must have $n_1<n$ and $V_1=V\cap
H^0(E_1)$ (otherwise replacing $V_1$ by $V\cap H^0(E_1)$ would
contradict the $\a_i$-semistability of $(E,V)$).

\n{\em Proof of (i).}  Suppose first that $(E,V)$ is
$\a_i^+$-stable but $\a^-_i$-unstable. We observe that if
$(E_1,V_1)\in {\cF}_1$, then $\frac{k_1}{n_1}< \frac{k}{n}$, since
otherwise $(E,V)$ could not be $\a_i^+$-stable. But the allowed
values for $\frac{k_1}{n_1}$ are limited by the constraints
$0<n_1< n$ and $0\le k_1\le k$. We can thus define
 $$
  \l_0=\max\left\{\frac{k_1}{n_1}\ \biggm|\ (E_1,V_1)\in {\cF}_1\ \right\}\
 $$
and set
 $$
 {\cF}_2=\left\{(E_1,V_1)\subset {\cF}_1\ \biggm|\\
 \frac{k_1}{n_1}=\l_0\ \right\}.
 $$

\n Let $(E_1,V_1)$ be any coherent system in $\cF_2$. Since
$V_1=V\cap H^0(E_1)$, we can write
 $$
 0\to (E_1,V_1)\to (E,V)\to (E_2,V_2) \to 0 \
 $$
for some coherent system $(E_2,V_2)$. Since  $\mu_{\a_i}(E_1,V_1)=
\mu_{\a_i}(E,V)=\mu_{\a_i}(E_2,V_2)$ and $(E,V)$ is
$\a_i$-semistable, it follows that both
$(E_1,V_1)$ and $(E_2,V_2)$ are $\a_i$-semistable.

We now show that $(E_2,V_2)$ is $\a_i^+$-stable. Suppose not. Then
there is a proper subsystem $(E_2',V_2')\subset (E_2,V_2)$ with
\begin{itemize}
\item $\mu_{\a_i}(E_2',V_2')=\mu_{\a_i}(E_2,V_2)$,
\item $\frac{k'_2}{n'_2}\ge\frac{k_2}{n_2}$.
\end{itemize}
\n Consider now the subsystem $(E',V')\subset (E,V)$ defined by
the pull-back diagram
 $$
 0\to (E_1,V_1)\to (E',V')\to (E'_2,V'_2) \to 0 \ .
 $$
This has $\mu_{\a_i}(E',V')=\mu_{\a_i}(E,V)$ and thus satisfies
$\frac{k'_2+k_1}{n'_2+n_1}\le\frac{k_1}{n_1}$. It follows that
 $$
 \frac{k'_2}{n'_2}\le\frac{k_1}{n_1}<\frac{k_2}{n_2}\ ,
 $$
which is a contradiction.

 Now consider $(E_1,V_1)\in\cF_2$ with {\it minimum rank in}
$\cF_2$. If $(E_1,V_1)$ is not $\a_i^+$-stable, then it must have a
proper subsystem $(E'_1,V'_1)$ with
\begin{itemize}
\item $\mu_{\a_i}(E_1',V'_1)=\mu_{\a_i}(E_1,V_1)$,
\item $\frac{k'_1}{n'_1}\ge\frac{k_1}{n_1}$.
\end{itemize}
\n But then $n'_1<n_1$, which contradicts the minimality of
$n_1$.  Finally, notice that since $(E,V)$ is $\a_i^+$-stable,
we must have
$\mu_{\a_i^+}(E_1,V_1)<\mu_{\a_i^+}(E,V)<\mu_{\a_i^+}(E_2,V_2)$.

\n{\em Proof of (ii).}   If $(E,V)$ is $\a_i^-$-stable but
$\a^+_i$-unstable, then $\frac{k_1}{n_1}> \frac{k}{n}$ for all
$(E_1,V_1)\in\cF_1$. The proof of (i) must thus be modified as
follows. With
 $$
 \l_0=\min\left\{\frac{k_1}{n_1}\ \biggm|\ (E_1,V_1)\in {\cF}_1\ \right\}\
 $$
we can define

 $$
 {\cF}_2=\left\{(E_1,V_1)\subset {\cF}_1\ \biggm|\\
 \frac{k_1}{n_1}=\l_0\ \right\}
 $$
and select $(E_1,V_1)\in\cF_2$ such that $E_1$ has minimal rank in
$\cF_2$. It follows in a similar fashion to that above that
$(E,V)$ has a description as
 $$
 0\to (E_1,V_1)\to (E,V)\to (E_2,V_2) \to 0
 $$
in which both $(E_1,V_1)$ and $(E_2,V_2)$ are $\a_i^-$-stable.
\qed

We refer to the extensions of the form (\ref{destab})
with the properties of Lemma \ref{lem:vicente} as the
{\em destabilising patterns} of the coherent systems.

%%%%%%%%%%%%%%%%%
\subsection{Codimension estimates for $G_i^-$ and $G_i^+$}
%%%%%%%%%%%%%%%%%

\begin{definition}\label{def:wi}\begin{em} Let $W^+(\a_i, \l, n_1; n,d,k)$
(abbreviated to $W^+_i(\l,n_1)$ whenever possible) denote the set
of all destabilising patterns
 $$
 0\to (E_1,V_1)\to (E,V)\to (E_2,V_2) \to 0
 $$
in which
\begin{itemize}
\item $(E,V)$ is  $\a^+_i$-stable and of type $\ndk$,
\item ${\rk}(E_1)=n_1$ and $\dim(V_1)=\l n_1$,
\item
$\mu_{\a_i}(E_1,V_1)=\mu_{\a_i}(E_2,V_2)=\mu_{\a_i}(E,V)$,
\item $(E_1,V_1)$ and $(E_2,V_2)$ are both $\a_i^+$-stable,
\item $\dim(V_1)$ and ${\rk}(E_1)$ satisfy the min-max criteria given in
{\rm (c)} and {\rm (d)} of Lemma \ref{lem:vicente}(i).
\end{itemize}

Similarly, let $W^-(\a_i, \l, n_1; n,d,k)$ (abbreviated to
$W^-_i(\l,n_1)$ whenever possible) denote the set of all
destabilising patterns
 $$
 0\to (E_1,V_1)\to (E,V)\to (E_2,V_2) \to 0
 $$
in which
\begin{itemize}
\item $(E,V)$ is  $\a^-_i$-stable and of type $\ndk$,
\item ${\rk}(E_1)=n_1$ and $\dim(V_1)=\l n_1$,
\item
$\mu_{\a_i}(E_1,V_1)=\mu_{\a_i}(E_2,V_2)=\mu_{\a_i}(E,V)$,
\item $(E_1,V_1)$ and $(E_2,V_2)$ are both $\a_i^-$-stable,
\item $\dim(V_1)$ and ${\rk}(E_1)$ satisfy the min-min criteria given in
{\rm (c)} and {\rm (d)} of Lemma \ref{lem:vicente}(ii).
\end{itemize}

Define
 $$
 W^+(\a_i,n,d,k)=\bigsqcup_{\l<\frac{k}{n},\,
 n_1<n} W^+(\a_i, \l, n_1; n,d,k),
 $$
 $$
 W^-(\a_i,n,d,k)=\bigsqcup_{\l>\frac{k}{n},\, n_1<n}
 W^-(\a_i, \l, n_1; n,d,k).
 $$
We abbreviate these to $W^+_i$ and $W^-_i$ whenever possible.\end{em}
\end{definition}

\begin{lemma}\label{lem:wi} Fix $\ndk$ and also $\a_i$. Then each set
$W^{\pm}_i(\l,n_1)$ is contained in a family of dimension bounded
above by

 \begin{eqnarray*} w^{\pm}_i(\l,n_1)&=
 \dim G(\a_i^{\pm};n_1,d_1,k_1)&\mbox{} +
 \dim G(\a_i^{\pm};n_2,d_2,k_2)\\
 & & \mbox{} + \max{\dim  \Ext^1((E_2,V_2),(E_1,V_1))}-1.
 \end{eqnarray*}
Here $n=n_1+n_2$, $d=d_1+d_2$ and $k=k_1+k_2$, and the maximum is
taken over all $(E_1,V_1)$, $(E_2,V_2)$ which satisfy the relevant
part of Definition \ref{def:wi}. Thus the set $W^{+}_i$ is
contained in a family whose dimension is bounded above by the
maximum of $w^{+}_i(\l,n_1)$  for all $\l<\frac{k}{n} $ and $\
n_1<n$. Similarly,  the set $W^{-}_i$ is contained in a family
whose dimension is bounded above by the maximum of
$w^{-}_i(\l,n_1)$  for all $\l>\frac{k}{n} $ and $\ n_1<n$.

\end{lemma}

\n{\em Proof.} In general the coherent systems moduli spaces
do not support universal objects. In order to obtain families in the
strict sense of the term, it is necessary to lift back from the
moduli spaces to a level
(for example, that of Quot schemes) on which families can be constructed.
One can then do a dimensional calculation. In fact this gives the same
answer is if we simply assumed that the moduli spaces support genuine
families (for a similar calculation, see, for example, \cite[Lemma 4.1]{BGN}).
Given this, the lemma follows at once from the definitions and Lemma
\ref{lem:vicente}.\hfill$\Box$

Note that $G(\a_i^+;n_1,d_1,k_1)=G_i(n_1,d_1,k_1)$ and
$G(\a_i^-;n_1,d_1,k_1)=G_{i-1}(n_1,d_1,k_1)$; the version used in
the lemma appears more natural in this context.

There are clearly surjective maps
 $$
 W^{\pm}_i\twoheadrightarrow G_i^{\pm}\ .
 $$
The maps may fail to be injective because a
coherent system in $G_i^{\pm}$ may have more than one subsystem which
satisfies the criteria on $(E_1,V_1)$ in Lemma \ref{lem:vicente}.
 Nevertheless, we can use the dimension
estimates on $W^{\pm}_i$ to estimate the codimension of $G_i^{\pm}$ in
$G(\a_i^{\pm};n,d,k)$ by
 $$
 \codim G_i^{+} \ge \dim
 G(\a^{+}_i;n,d,k) - \max\left\{ w^+_i(\l,n_1)\ \biggm|\
 \l<\frac{k}{n}\ ,\ \ n_1<n \right\}
 $$
 and
 $$
 \codim G_i^{-}
 \ge \dim G(\a^{-}_i;n,d,k) - \max\left\{ w^-_i(\l,n_1)\
 \biggm|\ \l>\frac{k}{n}\ ,\ \ n_1<n \right\}.
 $$

It follows from (\ref{destab}) and Proposition
\ref{prop:filtration}(ii) that in our situation
 $$
 \HH^0_{21}={\Hom}((E_2,V_2),(E_1,V_1))=0.
 $$
When $G(\a_i^{\pm};n,d,k)$, $G(\a_i^{\pm};n_1,d_1,k_1)$ and
$G(\a_i^{\pm};n_2,d_2,k_2)$ have their expected dimensions, and
$\mathbb H^2_{21}$ is zero for all relevant $(E_1,V_1)$ and
$(E_2,V_2)$, we have
\begin{eqnarray}
 \codim G_i^{+}&\geq&\b(n,d, k)\nonumber\\
 & &\mbox{}-\max\left\{ (\b(n_1,d_1,k_1) + \b(n_2,d_2,k_2) +
 C_{21}-1)\, \biggm|\,
 \frac{k_1}{n_1}<\frac{k}{n}\ ,\ \ n_1<n \right\}\nonumber \\
   &=&  \min\left\{\ C_{12}\ \biggm|\
   \frac{k_1}{n_1}<\frac{k}{n}\ ,\ \ n_1<n \right\}\label{codim+}
\end{eqnarray}
 by Corollary \ref{cor:euler}. Similarly
\begin{equation}\label{codim-}
\codim G_i^-
   \geq \min\left\{\ C_{12}\ \biggm|\
   \frac{k_1}{n_1}>\frac{k}{n}\ ,\ \ n_1<n \right\}.
\end{equation}

Of course in general we have to allow for the fact that the moduli
spaces may have dimensions greater than the expected ones and take
into account the contribution from  $\HH^2$ in the
computations of the actual dimensions. For later use, we state a
very general result and then we particularise to a result that
covers the cases considered in this paper.

In general, we shall describe the process of going from
$G(\a^+_i;n,d,k)$ to $G(\a^-_i;n,d,k)$ (or vice versa) as a {\em
flip}, although it is not necessarily a flip in any technical
sense. For all allowable values of $(\l,n_1)$, we denote the image
of $W^+_i(\l,n_1)$ in $G^+_i$ by $G^+_i(\l,n_1)$. For any
irreducible component $G$ of $G(\a^+_i;n,d,k)$, we shall say that
the flip is $(\l,n_1)$-{\em good on} $G$ if $G^+_i(\l,n_1)\cap G$
has positive codimension in $G$. A similar definition applies to
irreducible components of $G(\a^-_i;n,d,k)$. If a flip is
$(\l,n_1)$-good on all irreducible components of both
$G(\a^+_i;n,d,k)$ and $G(\a^-_i;n,d,k)$ and for all allowable
values of $(n_1,\l)$, we shall call it a {\em good flip}.

\begin{lemma}
\label{lem:flip-good}
  Let $\a_i$ be a critical value and suppose that
\begin{itemize}
\item $n_1+n_2=n$, $d_1+d_2=d$, $k_1+k_2=k$,
\item $\frac{d_1}{n_1}+\a_i\frac{k_1}{n_1}=
\frac{d_2}{n_2}+\a_i\frac{k_2}{n_2} =\frac{d}{n}+\a_i\frac{k}{n}$,
\item $\l=\frac{k_1}{n_1}<\frac{k}{n}$.
\end{itemize}

\n Let $G$ be an irreducible component of $G(\a^{+}_i;n,d,k)$ of
excess dimension $e\ge0$. Let $\{S_t\}$ be a stratification of
 $$
 G(\a^{+}_i;n_1,d_1,k_1)\x G(\a^{+}_i;n_2,d_2,k_2)
 $$
such that %$\dim\HH^0_{21}$ and $\dim\HH^2_{21}$ are both constant
  $\dim\HH^2_{21}$ is constant
on each $S_t$. Write $e_1$, $e_2$ for the excess dimensions of
irreducible components $G^1$ of $G(\a^{+}_i;n_1,d_1,k_1)$ and
$G^2$ of $G(\a^{+}_i;n_2,d_2,k_2)$. Then the flip at $\a_i$ is
$(\l,n_1)$-good if
 \begin{equation}\label{flipeq}
 C_{12}> %\dim\HH^0_{21}+
 \dim\HH^2_{21}+e_1+e_2-e-\codim_{G^1\x G^2}
 (S_t\cap(G^1\x G^2))
 \end{equation}
for all $G^1$, $G^2$ and all $S_t$ such that there exist
extensions (\ref{destab}) satisfying the conditions of Lemma
\ref{lem:vicente}(i) with $(E,V)\in G$ and
$((E_1,V_1),(E_2,V_2))\in S_t\cap(G^1\x G^2)$.

A similar result holds for $G(\a^{-}_i;n,d,k)$ if we replace the
condition $\l<\frac{k}{n}$ by $\l>\frac{k}{n}$ and Lemma
\ref{lem:vicente}(i) by Lemma \ref{lem:vicente}(ii).
\end{lemma}

\n{\em Proof.} We need to adjust the formulae (\ref{codim+}) and
(\ref{codim-}) by allowing for all the obstructions. For this we use
(\ref{C21}) and recall that we have already noted that
$\HH^0_{21}=0$.
\hfill$\Box$

\begin{corollary}\label{cor:flip}
Suppose that, for every allowable choice of $(n_1,d_1,k_1)$ with
$\frac{k_1}{n_1}<\frac{k}{n}$, $G(\a_i^{\pm};n_1,d_1,k_1)$ and
$G(\a_i^{\pm};n_2,d_2,k_2)$ have the expected dimensions, and that
stratifications $\{S_t^+\}$, $\{S_t^-\}$ of
$$G(\a^{+}_i;n_1,d_1,k_1)\x G(\a^{+}_i;n_2,d_2,k_2),\quad
 G(\a^{-}_i;n_1,d_1,k_1)\x G(\a^{-}_i;n_2,d_2,k_2)$$
exist such that $\dim\HH^2_{21}$
is constant on every stratum $S_t^+$ and $\dim\HH^2_{12}$ is constant on
every stratum $S_t^-$. Suppose further that
\begin{equation}\label{criteria}
C_{12}>\dim\HH^2_{21}-\codim S_t^+,\quad\mathrm{and}\quad
C_{21}>\dim\HH^2_{12}-\codim S_t^-
\end{equation}
for every $(n_1,d_1,k_1)$ and every stratum $S_t^{\pm}$. Then the flip
at $\a_i$ is good.
\end{corollary}

\n{\em Proof.} The hypotheses give $e_1=e_2=0$ for every
choice of $G^1$, $G^2$. %, while $\HH^0_{21}$ is always zero in our
%situation by Proposition \ref{prop:filtration}(ii) and (\ref{destab}).
The flip is therefore $(\l,n_1)$-good for $\l< \frac{k}{n}$ by
Lemma \ref{lem:flip-good}.

Now note that interchanging the indices $12$ changes a
destabilising pattern with $\l=\frac{k_1}{n_1}<\frac{k}{n}$ into
one with $\l=\frac{k_2}{n_2}>\frac{k}{n}$ and vice-versa. So the
second inequality in the statement shows that the flip is good for
$\l>\frac{k}{n}$.\hfill$\Box$

Of course, one needs to prove (\ref{criteria}) only for non-empty
strata. Moreover, if the extension (\ref{destab}) is trivial,
$(E,V)$ cannot be $\a$-stable for any $\a$. So, for proving the
first inequality, we may also assume that
$\dim\Ext^1((E_2,V_2),(E_1,V_1))>0$, i.e. by (\ref{C21})
 $$
 C_{21}+\dim\HH^2_{21}>0.
 $$
Similarly, for the second inequality, we may assume
 $$
 C_{12}+\dim\HH^2_{12}>0.
 $$

%%%%%%%%%%%%%%%%%%%%%%%%%%%%%%%%%%%%%%%%%%%%%%%%%%%%%%%%
\section{Coherent systems with $k=1$}\label{sec:k=1}
%%%%%%%%%%%%%%%%%%%%%%%%%%%%%%%%%%%%%%%%%%%%%%%%%%%%%%%%

We want to deal with applications of the theory developed so far to
the case of coherent systems with few sections and also to the case
of small rank. We devote the following sections to this task.

We start by analysing the case  $k=1$ and $n\geq 2$. The moduli space
of coherent systems in this case coincides with the moduli space
of pairs $(E,\p)$ which are $\a$-stable (see~\cite{BG1}). The
particular case $n=2$, $k=1$, $d>0$ has been studied thoroughly by
Thaddeus \cite{Th}, showing in particular that the spaces
$G(\a;2,d,1)$ are irreducible and of the expected dimension $2g
+d-2$. We assume that $g\ge2$ partly because of the complications of
Remark \ref{rmk:g=0,1} and partly because the proof fails for $g=0$.

\begin{theorem}[\cite{BD1,BD2,G,BDW,BDGW}]\label{thm:k=1}
  Let $g\ge2$. For $n>1$, the moduli spaces $G_i(n,d,1)$
  are non-empty, smooth, irreducible and of the expected dimension
  $\b=(n^2-n)(g-1)+d$. They are
  birationally equivalent for different values of $i$.
  The critical values are all of the form $\frac{s}{m}\in (0,\frac{d}{n-1})$
  with $0<m<n$ and $0<s<d$.
\end{theorem}

\n {\em Proof.\/}
  The smoothness property follows from Proposition \ref{prop:smooth}.
  Theorem \ref{G_L(k<n)} shows that the large $\a$ moduli space $G_L$
  is irreducible and of the
  expected dimension. So it only remains to prove that all the
  moduli spaces are birationally equivalent for different values
  of $\a$. This follows at once when we check that the flips are
  good. By Corollary \ref{cor:flip} we need only to verify the
inequalities (\ref{criteria}) for $k_1=0$, $k_2=1$, but we do need to know
that all non-empty $G(\a_i^{\pm};n_1,d_1,k_1)$ with $n_1<n$ and $k_1=0,1$
have the expected dimensions. For $k_1=0$, these spaces are the full
moduli spaces, for which we know the result to be true. We can therefore
proceed by induction on $n$.

For the base case, we take the equivalent theorem for $n=1$, namely that
$G(1,d,1)$ has dimension $d$. This is clear since $G(1,d,1)=S^dX$.

We can therefore proceed to the inductive step.   Note first that
$\HH^2_{21}=0$ by Lemma \ref{lem:hopf} and $\HH^2_{12}=0$ since
$V_1=0$. The critical value $\a_i$ is given by
 $$
 \frac{d_1}{n_1}=\frac dn +\frac{\a_i}{n}
    =\frac{d_2}{n_2}+\frac{\a_i}{n_2},
 $$
i.e.
\begin{equation}\label{eqn:critical}
\a_i=\frac{1}{n_1}(d_1n_2-d_2n_1)
 \end{equation}
We have by (\ref{dim-ext})
  $$
  C_{12}=n_1n_2(g-1)-d_2n_1+d_1n_2=n_1n_2(g-1)+n_1\a_i>0.
  $$
On the other hand
 $$
 C_{21}=n_1n_2(g-1)-d_1n_2+d_2n_1+d_1-n_1(g-1).
 $$
Now
 $$
 d_1n_2-d_2n_1=n_1\a_i<\frac{n_1d}{n-1},
 $$
which gives $d_1n_2-d_2n_1<d_1$. So
 $$
 C_{21}>n_1(n_2-1)(g-1)\ge0.
 $$
\hfill $\Box$
%%%%%%%%%%%%%%%%%%%%%%%%%%%%%%%%%%%%%%%%%%%%%%%%%%%%%%%%
\section{Coherent systems with $k=2$}\label{sec:k=2}
%%%%%%%%%%%%%%%%%%%%%%%%%%%%%%%%%%%%%%%%%%%%%%%%%%%%%%%%

We look next at the case $k=2$.
\begin{theorem} \label{thm:k=2}
  Let $X$ be a Petri curve of genus $g\geq 2$. Then we have
  \begin{itemize}
   \item  For $n=2$ the moduli spaces $G_i(2,d,2)$ are
    non-empty if and only if $d> 2$.
    They are irreducible and of the expected dimension $2d-3$.
   \item  For $n>2$ the moduli spaces $G_i(n,d,2)$ are
    non-empty if and only if $d>0$. They are always irreducible
    and of the expected dimension $(n^2-2n)(g-1)+2d-3$.
  \end{itemize}
\end{theorem}

\n {\em Proof.\/}
  We start by considering the moduli space $G_L$. Here the result follows
from
  Theorem \ref{G_L(k=n)} when $n=2$ and from
  Theorem \ref{G_L(k<n)} when $n>2$.

It remains to prove that all the flips are good. Again we proceed by
induction
on $n$, noting that we already know that the moduli spaces for $k=0,1$ do
have
the expected dimensions. For the base case, we take the statement that the
moduli
spaces $G(1,d,2)$ have the expected dimensions. This is true by section
\ref{subsec:alpha-small} since we are assuming that the curve is Petri.
Note incidentally that these spaces are not necessarily irreducible,
but irreducibility is not needed for the argument.

We now proceed to the inductive step. According to Corollary \ref{cor:flip},
we can restrict attention to the two cases $k_1=0$, $k_2=2$  and
$k_1=k_2=1$,
$n_1>n_2$. In each case we need to prove the inequalities (\ref{criteria}).

\begin{enumerate}
\item $k_1=0$, $k_2=2$. The critical value is given by
 $$
 \frac{d_1}{n_1}=\frac{d_2}{n_2}+\frac{2\a_i}{n_2},
 $$
i.e.
 $$
\a_i=\frac{1}{2n_1}(d_1n_2-d_2n_1).
 $$
 So
 $$
 C_{12}=n_1n_2(g-1)-d_2n_1+d_1n_2=n_1n_2(g-1)+2n_1\a_i>0.
 $$
  By Proposition \ref{prop:C21},
  $\HH^2_{21}=\Ext^2((E_2,V_2),(E_1,0))=H^0(E_1^*\ox N_2 \ox K)^*$,
  where $N_2$ is the kernel of ${\cO}^2\to E_2$. If $N_2=0$ we have
  finished as we have already proved that $C_{12}>0$. When $N_2$ is
  non-zero we have an exact sequence $N_2\to
  {\cO}^2\twoheadrightarrow L$ onto some line bundle $L$ with at
  least two sections. Therefore $\deg N_2=-\deg L \leq -\frac{g+2}2$,
  by section \ref{subsec:basicBN}, since the curve is Petri. So
  $$
  \deg (E_1^*\ox N_2\ox K) \leq -d_1+n_1(2g-2-\frac{g+2}2)
  <n_1(2g-2).
  $$
  Then by Clifford's theorem \cite{BGN} applied to the semistable
  bundle $E_1^*\ox  N_2\ox K$, if $h^0(E_1^*\ox  N_2\ox K)>0$ then
 \begin{eqnarray*}
  \dim \HH^2_{21} &\leq & \frac{-d_1+n_1(2g-2-\frac{g+2}2)}2 +n_1=
  -\frac{d_1}{2} + \frac{n_1}{4}(3g-2) \\
   & <& \frac34 n_1(g-2) + n_1 < n_1n_2(g-2)+ n_1n_2+
   2n_1\a_i=C_{12}.
 \end{eqnarray*}

  On the other hand, $\HH^2_{12}=0$ since $k_1=0$. Therefore
  we only need to prove that $C_{21}>0$. Now
 $$
  C_{21}=n_1n_2(g-1)+n_1 d_2-n_2d_1-2n_1(g-1)+2d_1.
 $$
  If $n_2>2$ then we use the bound on the
  $\a$-range given by $\a_i<\frac{d_2}{n_2-2}$. Hence $\a_i<
  \frac{d_2}{n_2}+\frac{2}{n_2}\a_i=\frac{d_1}{n_1}$ and
  $$
  C_{21}=n_1(n_2-2)(g-1)+2d_1-2n_1\a_i>0.
  $$
  If $n_2=2$ then
  $d_2> 2$ by induction hypothesis, and so $\frac{d_1}{n_1}=
  \frac{d_2}{2}+\a_i > \a_i$, whence $C_{21}=2d_1-2n_1\a_i >0$.
  If $n_2=1$ then $d_2\geq \frac{g+2}2$ since $E_2$ is a line
  bundle with at least two sections on a Petri curve. As
  $\frac{d_1}{n_1}>d_2$, we have
 \begin{eqnarray*}
  C_{21} &=& -n_1(g-1)+n_1d_2+d_1 > 2n_1d_2-n_1(g-1)\\
  &\geq & n_1(g+2-g+1)>0.
 \end{eqnarray*}
  So in all the cases $C_{21}>0$, as required.

\item $k_1=k_2=1$, $n_1>n_2$. The critical value is given by
  $$
  \frac{d_1}{n_1}+\frac{\a_i}{n_1}=\frac{d_2}{n_2}+\frac{\a_i}{n_2}=
  \frac{d}{n}+\frac{2}{n}\a_i.
  $$
  i.e.
  $$
  \a_i=\frac{1}{n_1-n_2}(d_1n_2-n_1d_2).
  $$
  By Lemma \ref{lem:hopf}, we have $\HH^2_{21}=0$ and
  $\HH^2_{12}=0$. We compute
 \begin{eqnarray*}
  C_{12}&=&n_1n_2(g-1)-n_1d_2+n_2d_1-n_2(g-1)+d_2-1=\\
  &=&(n_1-1)n_2(g-1)+\a_i(n_1-n_2)+d_2-1>0, \\
  C_{21} &=& n_1n_2(g-1)+n_1d_2-n_2d_1 -n_1(g-1)+d_1-1= \\
  &=& (n_2-1)n_1(g-1)+ d_1-\a_i(n_1-n_2)-1.
 \end{eqnarray*}
  For $n_2>1$ we use the $\a$-range condition to get
  $\a_i<\frac{d_1}{n_1-1}$ and so $d_1-\a_i(n_1-n_2)>d_1 -
  \frac{d_1}{n_1-1}(n_1-n_2)\geq 0$ and thus $C_{21}>0$. In the case
  $n_2=1$, we have $C_{21}=n_1d_2-1>0$.
\end{enumerate}
\hfill $\Box$

\begin{remark}\begin{em}\label{rem:G(2,2,2)}
 In the case $n=d=k=2$, $\tilde{G}_L$ consists only of reducible
 coherent systems and it is irreducible and of dimension $2$ by
 Theorem \ref{G_L(k=n)}. It
 is easy to see that in this case there are no flips.
\end{em}\end{remark}

%%%%%%%%%%%%%%%%%%%%%%%%%%%%%%%%%%%%%%%%%%%%%%%%%%
\section{Coherent systems with $n=2$}\label{sec:n=2}
%%%%%%%%%%%%%%%%%%%%%%%%%%%%%%%%%%%%%%%%%%%%%%%%%%

Now we are going to deal with coherent systems of rank $2$. Our
results in this case are partial. This is due to two reasons. On
the one hand our understanding of the moduli space $G_L$ of
coherent systems for large values of the parameter $\a$ for $k\geq
4$ is very limited, in particular we do not know whether these
spaces are irreducible and of the expected dimension. On the other
hand we only manage to check that the flips are good for $k\leq
4$. We need a preliminary result on rank $1$ coherent systems.

\begin{lemma}\label{lem:St}
  Let $X$ be a Petri curve of genus $g\geq 2$. Consider in $G(1,d,k)$ the
  stratification given by the sets $S_t=\{ (L,V) \in G(1,d,k)\;
  |\; h^0(L)=t\}$. Then
 \begin{itemize}
 \item If $d \leq g-1+k$ then the number of sections $h^0(L)$
  of a generic $(L,V) \in G(1,d,k)$ is $k$, and $\codim
  S_{k+j}=j(g-d-1+k+j)$, when non-empty.
 \item If $d \geq g-1+k$ then the number of sections $h^0(L)$
  of a generic $(L,V) \in G(1,d,k)$ is $p=d-g+1$, and $\codim
  S_{p+j}=j(d-g+1-k+j)$, when non-empty.
 \end{itemize}
\end{lemma}

\n {\em Proof. \/}
  Let $p$ be the number of sections of a generic $(L,V) \in G(1,d,k)$. Then
  it must be $\dim G(1,d,k)=\dim G(1,d,p)+\dim \Gr(k,p)$. By an easy
computation
  it follows that either $p=d-g+1$ or $p=k$. If $d<g-1+k$ then it must be
  $p=k$ and $\codim S_{k+j}=\dim G(1,d,k)-\dim G(1,d,k+j)- \dim \Gr(k,k+j)$.
  If $d \geq g-1+k$ then $\codim S_{d-g+1}=0$ so $p=d-g+1$. The computation
  of $\codim S_{p+j}$ is left to the reader.
\hfill $\Box$

Now we focus on the study of $G_i(2,d,k)$ for $k>0$. The expected
dimension is $\b(2,d,k)=(4-2k)g + k d -k^2+2k- 3$. For $k=1$ this
has been treated in section \ref{sec:k=1} and for $k=2$ in section
\ref{sec:k=2}. So we may restrict to the case $k>2$. By Lemma
\ref{lem:degree} it must be $d>0$ for stable objects to exist.

\begin{theorem} \label{thm:n=2}
  Let $X$ be a Petri curve of genus $g\geq 2$. Then
  \begin{itemize}
  \item For $k=2$ the moduli spaces $G_i(2,d,2)$ are non-empty
  if and only if $d> 2$. They are irreducible
  and of the expected dimension $\b=2d-3$.
  \item For $k=3$ the moduli spaces $G_i(2,d,3)$ are non-empty
  if and only if $d\geq \frac{2g+6}3$. They are always of the
  expected dimension $\b=3d-2g-6$ and irreducible when $\b>0$.
  \item For $k=4$ the moduli spaces $G_i(2,d,4)$
  are birational to each other.
  \end{itemize}
\end{theorem}

\n {\em Proof.\/}
 We start by considering the moduli space $G_L$. Here the
 result follows from Theorem \ref{G_L(k=n)} for $k=2$ and from
 Theorem \ref{thm:dual-span} for $k=3$.

 Let now $k=2$, $3$ or $4$ and we will prove that the flips
 are good. By Corollary \ref{cor:flip} we have to prove the
 inequalities \eqref{criteria} for $n_1=n_2=1$ and all possible
 choices of $k_1<\frac{k}{2}$, since the moduli spaces of
 coherent systems of type $(1,d',k')$ have the expected
 dimension for a Petri curve, by section \ref{subsec:alpha-small}.
 As $k\leq 4$ we have that $k_1=0$ or $1$.

 More in general, let $k\geq 2$ be an integer, and consider
 extensions as in \eqref{destab} of the form
 $(L_1,V_1) \to (E,V) \to (L_2,V_2)$ where
 $n_1=n_2=1$ and $k_1<\frac{k}{2}$ satisfying $k_1\leq 1$.
 Then we are going to prove that the inequalities \eqref{criteria}
 are satisfied. By Lemma \ref{lem:flip-good} this implies that
 the flip is $(\l, 1)$-good on $G(\a^+_i;2,d,k)$ for
 $\l=0,1$ and $(\l,1)$-good on $G(\a_i^-;2,d,k)$ for
 $\l=k,k-1$.

 The critical value $\a_i$ is given by
 $$
  d_1+k_1\, \a_i=\frac d2+\frac k2\a_i= d_2 +k_2\, \a_i,
 $$
  i.e.
 $$
 \a_i=\frac{d_1-d_2}{k_2-k_1}.
 $$

 We start by proving the second inequality in \eqref{criteria}.
 In this case Lemma \ref{lem:hopf} implies that $\HH_{12}^2=0$
 since $k_1\leq 1$.
 By Theorem \ref{thm:BN}, in order for coherent systems
 of type $(1,d_2,k_2)$ to exist we must have
 \begin{equation}\label{eqn:T}
   d_1>d_2 \geq \frac{k_2-1}{k_2}g +k_2-1.
 \end{equation}
  We compute
 \begin{eqnarray*}
  C_{21}&=&g-1+d_2-d_1+k_2(d_1-g+1-k_1) \\
  &=&d_2 +(k_2-1)(d_1-g+1-k_1)-k_1 \\
  &\geq& k_2\left(\frac{k_2-1}{k_2}g +k_2-1\right)
  +(k_2-1)(-g+2-k_1)-k_1\\
  &=& (k_2-k_1-1)k_2+2(k_2-1) \geq 2k_1 > 0.
\end{eqnarray*}

 Now we prove the first inequality in \eqref{criteria}.
 We have
\begin{eqnarray*}
  C_{12}&=&g-1+d_1-d_2+k_1(d_2-g+1-k_2)\\
  &\geq& g-1+1 +k_1\left( \frac{k_2-1}{k_2}g
  -g\right)=\frac{k_2-k_1}{k_2}g>0.
\end{eqnarray*}
 If $\HH^2_{21}=0$ then we have finished.
 Otherwise, Lemma \ref{lem:hopf} gives the bound
 $$
  \dim \HH^2_{21}\leq (k_2-1)(h^0(L_1^*\ox K)-1).
 $$
 We stratify $G(1,d_1,k_1)$ by using the subsets defined in Lemma
 \ref{lem:St}. Let $S_t$ be the subspace of those $(L_1,V_1)\in
 G(1,d_1,k_1)$ with $h^0(L_1)=t$. It only remains to check that
 $C_{12}>\dim \HH^2_{21}-\codim S_t$ at
 the points in $S_t$.

 Suppose first that $d_1> g-1+k_1$. Lemma \ref{lem:St} says that
 the generic number of sections $h^0(L_1)$ of an element $(L_1,V_1)\in
 G(1,d_1,k_1)$ is $p=d_1-g+1$ and that $\codim S_{p+t}=t(d_1-g+1-k_1+t)$.
 Also $h^0(L_1^*\ox K)=t$ at a point in $S_{p+t}$.
 Suppose that $\dim \HH^2_{21}-\codim S_{p+t}>0$ since otherwise there
 is nothing to prove. So
\begin{eqnarray}\label{eqn:n=2.case1}
 \dim {\HH}^2_{12}-\codim S_{p+t} &\leq& (k_2-1)(t-1)-
 t(d_1-g+1-k_1+t) \nonumber \\
 &\leq & t(g-d_1-2+k-t)  \\
 &\leq& (k_2-1)(g-d_1+k-3), \nonumber
\end{eqnarray}
 since it must be $1\leq t \leq k_2-1$ for the
 second line to be non-negative. In the other case, $d_1\leq g-1+k_1$,
 the generic number of sections of $L_1$ is $p=k_1$ and $\codim
 S_{p+t}=t(g-1-d_1+k_1+t)$. Since $h^0(L_1^*\ox K)=g-1-d_1+k_1+t$
 at a point in $S_{p+t}$, we have
\begin{eqnarray}\label{eqn:n=2.case2}
 \dim {\HH}^2_{12}-\codim S_{p+t} &\leq& (k_2-1)(g-1-d_1+k_1+t-1)-
  t(g-1-d_1+k_1+t) \nonumber \\
 &\leq& (k_2-1-t)(g-1-d_1+k_1+t) \\
 &\leq&  (k_2-1)(g-d_1+k-3). \nonumber
\end{eqnarray}
 So using either \eqref{eqn:n=2.case1} or
 \eqref{eqn:n=2.case2} it only remains to prove that
  $$
 g-1+d_1-d_2+k_1(d_2-g+1-k_1k_2)>(k_2-1)(g-d_1+k-3).
  $$
 Rearranging terms this is equivalent to
  $$
 k_2d_1+(k_1-1)d_2>(k-2)g+(k-2)(k_2-2)+k_1k_2.
  $$
 Using \eqref{eqn:T} it suffices to show that
  $$
 k_2+(k-1)\left(\frac{k_2-1}{k_2}g
 +k_2-1\right)>(k-2)g+(k-2)(k_2-2)+k_1k_2.
  $$
 This holds for $k_1=0$ or $1$ and $k_2=k-k_1$.
\hfill $\Box$

\begin{remark} \begin{em}\label{rem:noflips}
  In order to have any flips, \eqref{eqn:T} imposes the condition
  $$
   d\geq 2\left( \frac{k_2-1}{k_2}g +k_2-1\right) +1,
  $$
  for some $k_2>\frac k2$. This implies that $d \geq 2(\frac g2 +1)+1=g+3$.
  So when $d\leq g+2$ there are no flips for $G(\a;2,d,k)$.
\end{em}\end{remark}

Checking whether the flips are good when $k_1>1$ is difficult in
general. Nonetheless we have the following positive result for the
case $k_1=2$.

\begin{theorem} \label{thm:k=5,6}
  Let $X$ be a Petri curve of genus $g\geq 2$. Consider the
  moduli spaces of coherent systems of type $(2,d,k)$ with $k>4$,
  and let $\a_i$ be a critical value corresponding to coherent
  subsystems with $n_1=1$ and $k_1=2$. Then the flip at $\a_i$ is
  $(\l=k-2,n_1=1)$-good on $G(\a_i^-;2,d,k)$.
  In particular, when $k=5$ or $k=6$, if the moduli space $G_0(2,d,k)$
  is non-empty then $G_L(2,d,k)$ is non-empty also.
\end{theorem}

\n {\em Proof.\/}
 By Lemma \ref{lem:flip-good} we need to check that for
 $n_1=n_2=1$ and $k_1=2$ we have the inequality
 $C_{21}>\dim \HH^2_{12}-\codim S_t$ at the points of $S_t$,
 for a suitable stratification $\{S_t\}$ of
 $G(\a^{-}_i;1,d_1,2)\x G(\a^{-}_i;1,d_2,k-2)$.
 By the proof of Theorem \ref{thm:n=2}, we
 already know that $C_{21}\geq 2k_1=4>0$.

 We distinguish two cases. First suppose that $d_2 \geq g+k-3$. We
 consider the stratification of $G(1,d_2,k-2)$ given by $S_t=\{
 (L_2,V_2)\; |\; h^0(L_2)=t\}$. By Lemma \ref{lem:hopf} we know that
 $\dim \HH^2_{12} \leq h^0(L_2^* \ox K)-1$. The
 generic number of sections of $L_2$ for an element $(L_2,V_2)
 \in G(1,d_2,k-2)$ is
 $p=d_2-g+1$. Using Lemma \ref{lem:St} we have that for any
 $t\geq 0$, at a point in $S_{p+t}$,
  $$
   \dim \HH^2_{12} -\codim S_{p+t} \leq t-1 -t(d_2-g+1-k+2+t) \leq
   0<C_{21},
  $$
 as required.

 The other case is $d_2 <g+k-3$. Then the generic number of
 sections of $L_2$ for an element $(L_2,V_2)$ is $p=k-2$.
 So for any $t\geq 1$  we have at a point in $S_{p+t}$,
 $$
 \dim \HH^2_{12} -\codim S_{p+t} \leq (1-t)(g-1-d_2+k-2+t) \leq 0<C_{21}.
 $$
 This means that we may restrict to the case where
 $(L_2,V_2)$ lies in the open subset $S_{k-2}\subset
 G(1,d_2,k-2)$. A coherent system $(L_2,V_2)\in S_{k-2}$ is
 determined by its underlying line bundle $L_2$. Now consider the
 exact sequence $N_1 \to {\cO}^2\to L_1$, where $N_1$ is the
 kernel. Then $N_1$ is a line bundle of degree $-l$, say.
 One clearly has $l\leq d_1$. Define
 the stratification of $S_{k-2}$ given by the subsets
 $$
  T_t=\{ (L_2,V_2) \in S_{k-2} \; |\; h^0(N^*_1\ox L_2)=t\}.
 $$
 Clearly $\dim T_t \leq \dim G(1,d_2+l,t)$. Also we stratify
 $G(1,d_1,2)$ by the subsets $W_l$ of those coherent
 systems $(L_1,V_1)$ such that the image of the map ${\cO}^2\to L_1$
 is a line bundle of degree $l$. Generically this map is
 surjective, so $W_{d_1}$ is an open dense subset.

 We start by considering the stratum $W_{d_1}\subset G(1,d_1,2)$.
 An easy calculation using that $d_1>d_2\geq \frac{k-3}{k-2}g +k-3$
 (see \eqref{eqn:T}) and $k\geq 5$
 shows that
  $$
  \dim G(1,d,d-g+4)< \dim G(1,d_2,k-2).
  $$
 Therefore the generic number of sections of the line bundle
 $N^*_1\ox L_2$, for $(L_2,V_2) \in S_{k-2}$ and $(L_1,V_1) \in W_{d_1}$,
 is $p\leq d-g+3$. Note that in particular $d-g+3 \geq 0$.
 At a point of $T_t\subset S_{k-2}$ with
 $t\leq d-g+3$ we have
  $$
   \dim \HH^2_{12}=\dim H^0(L_2^*\ox N_1\ox K)=g-1-d+t \leq 2<
   C_{21}.
  $$
 At a point of $T_{d-g+4+t}$ with $t\geq 0$ we have,
 \begin{eqnarray*}
  \dim \HH^2_{12} &- & \codim T_{d-g+4+t}=g-1-d+d-g+4+t -
  \codim T_{d-g+4+t} \\ &\leq & 3+t- \dim G(1,d_2,k-2)
  +\dim G(1,d,d-g+4+t) \\ &<& 3+t-t(d-g+7+t) \leq 3 \leq C_{21}.
 \end{eqnarray*}

 For the stratum $W_{d_1-1}$ we use that $\dim G(1,d-1,d-g+3)< \dim
 G(1,d_2,k-2)$ to prove that the generic number of sections of
 $N^*_1\ox L_2$, for $(L_2,V_2) \in S_{k-2}$ and
 $(L_1,V_1) \in W_{d_1-1}$, is $p\leq d-g+2$.
 Working as before we get that
  $$
  \dim \HH^2_{12}- \codim T_{d-g+2+t} \leq 2<C_{21},
  $$
 for $t\geq 0$.

 Finally consider the strata $W_l\subset G(1,d_1,2)$ where $l\leq
 d_1-2$. It is easy to check that $\codim W_l=d_1-l \geq 2$.
 We have that $l\geq \frac{g+2}2$, since $N^*_1$ has two
 sections and $X$ is a Petri curve. Now an easy calculation
 shows that
  $$
   \dim G(1,d_2+l, d_2+l-g+7)<\dim G(1,d_2,k-2).
  $$
 Therefore the generic number of sections of $N^*_1\ox L_2$ is $p\leq
 d_2+l-g+6$. At a point $((L_1,V_1),(L_2,V_2))
 \in W_l\x T_t \subset G(1,d_1,2) \x S_{k-2}$ with
 $t\leq d_2+l-g+6$ we have
 $$
  \dim \HH^2_{12} -\codim W_l \leq g-1-(d_2+l) +t -2 \leq 3<C_{21}.
 $$
 At a point of $W_l\x T_{d_2+l-g+7+t}$ with $t\geq 0$, we have
 $$
  \dim \HH^2_{12}-\codim W_l-\codim T_{d_2+l-g+7+t}
  < 6+t-t(d_2+l-g+13+t)-2 \leq 4\leq C_{21},
 $$
 concluding that in all cases the flip is
 $(k-2,1)$-good. \hfill $\Box$

%%%%%%%%%%%%%%%%%%%%%%%%%%%%%%%%%%%%%%%%%%%%%%%%%%%%%%%%%
\section{Coherent systems with $k=3$}\label{sec:k=3}
%%%%%%%%%%%%%%%%%%%%%%%%%%%%%%%%%%%%%%%%%%%%%%%%%%%%%%%%%

Now we shall work out the case of the moduli spaces $G_i(n,d,3)$
of coherent systems with $k=3$ sections and rank $n>1$. Note that
the case $n=2$ follows from section \ref{sec:n=2}. We need a
preliminary result, similar in spirit to Lemma \ref{lem:St} but
for the case of bundles of higher rank. This result is somewhat
restricted as the only input is information on coherent systems
with at most $2$ sections.

\begin{lemma}\label{lem:bound}
  Let $d\leq n(g-1)$. Stratify the moduli space $M(n,d)$ by
  $S_t=\{F\in M(n,d)\; | \; h^0(F)=t\}$.
  Then $2h^0(F^*\ox K)-\codim S_t\leq 2(n(g-1)-d)+1$ at
  a point in $S_t$.
\end{lemma}

\n {\em Proof. \/}
 For $F\in S_0$ we have $2h^0(F^*\ox K)= 2(n(g-1)-d)$. For
 $t=1$ we have, by Proposition \ref{small-alpha},
 $\dim S_1\leq \dim G(\a;n,d,1)= (n^2-n)(g-1)+d$,
 where $\a>0$ is a small number. Hence $\codim S_1\geq
 n^2(g-1)+1- (n^2-n)(g-1)+d=n(g-1)-d+1$ and
  $$
 2h^0(F^*\ox K)-\codim S_1 \leq 2(n(g-1)-d+1)-n(g-1)-d+1=n(g-1)-d+1.
  $$
 For $t\geq 2$ we have that
 $\dim S_t +\dim \Gr(2,t) \leq \dim
 G(\a;n,d,2)=(n^2-2n)(g-1)+2d-3$, using
 Theorem \ref{thm:k=2}. So we deduce that
  \begin{eqnarray*}
  & &2h^0(F^*\ox K)- \codim S_t \leq  \\
  & & \leq 2(n(g-1)-d+t)-n^2(g-1)-1+(n^2-2n)(g-1)+2d-3-2(t-2)=0.
  \end{eqnarray*}
  The statement follows.
\hfill $\Box$

Now we obtain Clifford bounds type results for coherent systems.
The following results are not sharp, but they are good enough for
our purposes in this section. In the next two Lemmas, $X$ is {\em
any\/} curve of genus $g\geq 2$.

\begin{lemma} \label{clif-lem}
  Suppose $(E,V)$ is an $\a$-semistable coherent system with
  $\mu(E) \geq 2g-2$ and $h^1(E) >0$. Then
  $$
  h^0(E) \leq \frac{d}{2} +n + (n-1)k\a.
  $$
\end{lemma}

\n {\em Proof. \/}
  We want to bound $h^0(E)=h^1(E)+ d+n(1-g)$. Put
  $N=h^1(E)=h^0(E^*\ox K)$. Then there are $N$ linearly independent
  maps $E \to K$. For any divisor $D$ on $X$ of degree $[\frac{N-1}{n}]$
  we may find a non-zero map $E\to K(-D)$. The $\a$-semistability
  implies then
\begin{eqnarray*}
 & & \frac{d}{n}+\a \frac{k}{n} \leq 2g-2-\deg D +\a k, \\
 & & \left[\frac{N-1}{n}\right]\leq 2g-2-\frac{d}{n}+\a k\frac{n-1}{n}, \\
 & & N \leq n(2g-2)-d+\a k (n-1)+n, \\
 & & h^0(E) \leq n(g-1)+n +\a k (n-1) \leq \frac{d}{2}+n+\a
 k(n-1).
\end{eqnarray*}
\hfill $\Box$

\begin{lemma}\label{lem:clif-cs}
  Let $(E,V)$ be an $\a$-semistable coherent system with
  $0\leq\mu(E)<2g-2$. Then
  $$
  h^0(E) \leq \frac{d}{2} +n +(n-1)k\a.
  $$
\end{lemma}

\n {\em Proof. \/}
  For $n=1$ the last term is dropped and the result is the usual
  Clifford theorem for line bundles. Also for $\a>0$ very
  small, $E$ is a semistable bundle and the result follows by the
  Clifford theorem in \cite{BGN}. We also may suppose that $k>0$.
  Note that the bound weakens as we increase $\a$, so it is
  enough to check what happens when we cross a critical value
  $\a_i$ to the coherent systems $(E,V)$ that are
  $\a_i$-semistable but not $\a_i^-$-semistable.
  Then there is a pattern
 $$
  0\to (E_1,V_1) \to (E,V) \to (E_2, V_2) \to 0,
 $$
  with $k_1/n_1 < k/n < k_2/n_2$, $\mu_{\a_i}(E_1, V_1)=
  \mu_{\a_i}(E_2, V_2)=\mu_{\a_i}(E, V)$ and where $(E_1,V_1)$,
  $(E_2, V_2)$ are $\a_i$-semistable. Therefore $k_2>0$, and by
  Lemma \ref{lem:degree} we have $d_2 \geq 0$. Hence $0\leq
  d_2/n_2< d/n<2g-2$, and by induction,
 \begin{equation} \label{clif2}
    h^0(E_2) \leq \frac{d_2}{2} +n_2 +(n_2-1) k_2\a_i.
 \end{equation}
  There are three cases to consider:
\begin{itemize}
 \item $d_1/n_1 <2g-2$. As $d_1/n_1 > d/n \geq 0$, we apply
  induction to get
 $$
    h^0(E_1) \leq \frac{d_1}{2} +n_1 +(n_1-1)k_1\a_i\, ,
 $$
  which together with \eqref{clif2} gives the result using
  that $h^0(E)\leq h^0(E_1)+h^0(E_2)$. Note that $(n_1-1)k_1 +
  (n_2-1)k_2 \leq (n-1)k$, whenever $n=n_1+n_2$, $0<n_1,n_2<n$ and
  $k=k_1+k_2$, $k_1,k_2\geq 0$.
 \item $d_1/n_1 \geq 2g-2$ and $h^1(E_1) \neq 0$.
  We use Lemma \ref{clif-lem} to conclude
$$
    h^0(E_1) \leq \frac{d_1}{2} +n_1 +(n_1-1)k_1\a_i\, ,
$$
  and the result follows as in the previous case.
\item $d_1/n_1 \geq 2g-2$ and $h^1(E_1) =0$. Then
  $h^0(E_1)=d_1+n_1(1-g)$. We have
 $$
  h^0(E_1) =\frac{d_1}{2} + \frac{n_1}{2}\left[ \frac{d}{n}+
  \a_i\left(\frac{k}{n}-\frac{k_1}{n_1}\right)\right]+
  n_1(1-g) <
 $$
 $$
  <\frac{d_1}{2}+ \frac{n_1}{2}(2g-2) +n_1(1-g) +
  \a_i \frac{n_1k}{2n}< \frac{d_1}{2}+n_1+\a_i\frac{(n-1)k}{2n},
 $$
  from which we get again the result since $\frac{(n-1)k}{2n} +(n_2-1)k_2
  \leq (n-1)k$.
\end{itemize}
\hfill $\Box$

\begin{theorem} \label{thm:k=3}
  Let $X$ be a Petri curve of genus $g\geq 2$. Then we have
  \begin{itemize}
\item For $n=2$ the moduli spaces $G_i(2,d,3)$ are non-empty if and
  only if $d\geq \frac{2g+6}3$. They are always of the
  expected dimension $\b=3d-2g-6$ and irreducible when $\b>0$.
\item For $n=3$ the moduli spaces $G_i(3,d,3)$ are non-empty
  if and only if $d>3$.
  They are irreducible and of the expected dimension $\b=3d-8$.
\item For $n>3$ the moduli spaces $G_i(n,d,3)$ are non-empty if
  and only if $d>0$ and $d\geq n-(n-3)g$. They are always
  irreducible and of the expected dimension $\b=(n^2-3n)(g-1)+3d-8$.
  \end{itemize}
\end{theorem}

\n {\em Proof.\/}
  The case $n=2$ follows from Theorem \ref{thm:n=2}, so we may
  restrict to the case $n\geq 3$. The moduli space $G_L$ for the
  largest possible values of the parameter satisfies the
  statement of the Theorem, using Theorem \ref{G_L(k=n)} for the
  case $k=n=3$ and Theorem \ref{G_L(k<n)} for the case $n>k=3$.

  It remains to check that the flips are good. We proceed by
  induction on $n$, noting that we already know that the moduli
  spaces for $k=0,1,2$ have the expected dimensions for a Petri
  curve. According to Corollary \ref{cor:flip}, we have two cases:
  $k_1=0$, $k_2=3$ and $k_1=1$, $k_2=2$.

\begin{enumerate}
\item $k_1=0$, $k_2=3$. The critical value $\a_i$ is given by
 $$
 \frac{d_1}{n_1}=\frac{d_2}{n_2}+\frac{3}{n_2}\a_i=
 \frac{d}{n}+\frac{3}{n}\a_i\, .
 $$
 i.e.
 $$
 \a_i=\frac{d_1n_2-d_2n_1}{3n_1}\, .
 $$
 We start by proving the first inequality in \eqref{criteria}. We
 have
 $$
 C_{12}=n_1n_2(g-1)-n_1d_2+n_2d_1=n_1n_2(g-1)+3n_1\a_i>0.
 $$
 Now Lemma \ref{lem:hopf} implies $\dim \HH^2_{21} \leq
  2(h^0(E_1^*\ox K)-1)$ or else $\HH^2_{21}=0$. There are two cases:
\begin{enumerate}
 \item If $d_1 \leq n_1(g-1)$ then we use Lemma \ref{lem:bound}.
   Define the stratification given by
   $S_t=\{E_1\in M(n_1,d_1)\; |\; h^0(E_1)=t\}$. Then
   $$
    2h^0(E_1^*\ox K)-\codim S_t\leq 2(n_1(g-1)-d_1)+1.
   $$
   Hence $C_{12}> \dim \HH^2_{21}-\codim S_t$
   is implied by
   $$
    n_1(n_2-2)(g-1)+2d_1+3n_1\a_i>-1.
   $$
   For $n_2\geq 2$ this obviously holds.
   For $n_2=1$ we have that
    $$
    \qquad -n_1(g-1)+2d_1+3n_1\a_i=-n_1(g-1)+2d_1+d_1-n_1d_2=
    $$
    $$
    \qquad  =3d_1-n_1(g-1+d_2)\geq n_1(2d_2-g+1)>-1,
    $$
   using that $\frac{d_1}{n_1}>d_2\geq \frac{2g+6}3$, the last
   inequality being necessary for the existence of coherent
   systems of type $(1,d_2,3)$ on a Petri curve.
 \item If $d_1 > n_1(g-1)$ then we use Clifford theorem for the stable
   bundle $E_1^*\ox K$. So either $h^0(E_1^*\ox K)=0$ in which case
   there is nothing to prove, or
    $$
    h^0(E_1^*\ox K)\leq \frac{n_1(2g-2)-d_1}{2}+n_1,
    $$
   whence $\dim \HH^2_{21}\leq 2n_1g-d_1-2$. The inequality
   $C_{12} >2n_1g-d_1-2$ is equivalent to
    $$
    n_1(n_2-2)(g-2) + n_1(n_2-4) +d_1+3n_1\a_i >-2.
    $$
   For $n_2\geq 4$ this is obviously true. For $n_2=3$ it must
   be $d_2>3$ by induction hypothesis,
   so $\frac{d_1}{n_1} > 1+\a_i$ and $d_1 -n_1>0$, which
   yields the result.
   For $n_2=2$ we have $\frac{d_1}{n_1} > \frac{d_2}2 \geq 2$ as
   $d_2\geq \frac{2g+6}3$, by induction hypothesis.
   So $-2n_1+d_1 >0$ and we are done.
   For $n_2=1$ and $g\leq 5$ we have that $\frac{d_1}{n_1} >d_2\geq
   \frac{2g+6}3$ implies $\frac{d_1}{n_1}>d_2 \geq g+1$ and hence
    $$
     -n_1(g-2)-3n_1+d_1+3n_1\a_i> d_1-n_1(g+1)>-2,
    $$
   as required. The same argument covers the case $n_2=1$ and
   $d_1\geq n_1(g+1)$.
   Finally the case $n_2=1$, $n_1(g-1)<d_1<n_1(g+1)$ and $g\geq 6$
   requires a special
   treatment. We use the improvement of Clifford theorem given in
   \cite[Theorem 1]{M3}. Since $2+\frac{2}{g-4} \leq g-3
   <2g-2-\frac{d_1}{n_1}<g-1$ and the curve is Petri, we have
   $$
    h^0(E_1^*\ox K)\leq \frac{n_1(2g-2)-d_1}{2},
   $$
   which gives $\dim \HH^2_{21}\leq 2n_1(g-1)-d_1-2$, and hence
   $$
    C_{12}=n_1(g-1)+3n_1\a_i > \dim \HH^2_{21}.
   $$
\end{enumerate}

\n Now we pass on to prove the second inequality in
\eqref{criteria}. In this case $\HH^2_{12}=0$. We compute
 \begin{eqnarray*}
  C_{21}&=&n_1n_2(g-1)-3n_1 \a_i+3(d_1-n_1(g-1)) \\
  &=&n_1(n_2-3)(g-1)+3d_1-3n_1\a_i.
 \end{eqnarray*}
We have the following cases:
\begin{enumerate}
 \item If $n_2>3$ then $\a_i<\frac{d_2}{n_2-3}$. Computing we
 obtain that $\a_i<\frac{d_2}{n_2}+\frac{3}{n_2}\a_i=\frac{d_1}{n_1}$ and
 thus $C_{21}>0$.
 \item If $n_2=3$ then $C_{12}=3d_1-3n_1\a_i =d_2n_1 >0$.
 \item If $n_2=2$ then $\frac{d_1}{n_1}>\frac{d_2}2$. As
 $d_2\geq \frac{2g+6}{3}$ by induction hypothesis, we have
  \begin{eqnarray*}
   \qquad\qquad C_{21}&=&-n_1(g-1)+3d_1-2d_1+d_2n_1=n_1(d_2-g+1)+d_1 \\
   &>&n_1 \left(\frac32 d_2-g+1\right) \geq n_1(g+3-g+1)>0.
  \end{eqnarray*}
 \item If $n_2=1$ then  $d_2\geq \frac{2g+6}{3}$ in order to have
  stable coherent systems of type $(1,d_2,3)$ on a Petri curve. Also
  $\frac{d_1}{n_1}>d_2$, so
  \begin{eqnarray*}
   \qquad\qquad C_{21}&=&-2n_1(g-1)+3d_1-d_1+d_2n_1=n_1(d_2-2g+2)+2d_1 \\
   &>& n_1(3d_2-2g+2) \geq n_1(2g+6-2g+2)>0.
  \end{eqnarray*}
\end{enumerate}

\item $k_1=1$, $k_2=2$. The critical value is given by
  $$
  \frac{d_1}{n_1}+\frac{1}{n_1}\a_i=\frac{d_2}{n_2}+\frac{2}{n_2}\a_i=
  \frac{d}{n}+\frac{3}{n}\a_i\, ,
  $$
 i.e.
  $$
  \a_i=\frac{d_2n_1-d_1n_2}{n_2-2n_1}\, .
  $$
 It must be $n_2-2n_1 \neq 0$. We start proving the second inequality
 in \eqref{criteria}. We have $\HH^2_{12}=0$ and
  $$
  C_{21}=n_1n_2(g-1)+n_1d_2-n_2d_1+2(d_1-n_1(g-1)-1).
  $$
 We have the following cases:
\begin{enumerate}
 \item $n_2-2n_1>0$. Then $C_{21}=n_1(n_2-2)(g-1) +\a_i(n_2-2n_1)+2d_1-2>0$
  since $d_1>0$ and $n_2>2$.
 \item $n_2-2n_1<0$ and $n_2\geq 2$. Then
$\a_i=\frac{d_1n_2-d_2n_1}{2n_1-n_2}
  <\frac{d_1}{n_1-1}$ implies that $\a_i(2n_1-n_2)<2d_1-d_2$. So
  \begin{eqnarray*}
    \qquad C_{21} & =& n_1(n_2-2)(g-1)-\a_i(2n_1-n_2)+2d_1-2 \\
     &> & n_1(n_2-2)(g-1)+d_2-2 \geq 0.
  \end{eqnarray*}
  (Recall that in the particular case $n_2=2$ we must have $d_2> 2$.)
 \item $n_2-2n_1<0$ and $n_2=1$. Using that
  $\frac{d_1}{n_1}>d_2$ and $d_2\geq \frac{g+2}{2}$, we have
   \begin{eqnarray*}
    \qquad C_{21}&=& n_1(d_2-g+1)+d_1-2>n_1(2d_2-g+1)-2 \\
    &\geq & n_1(g+2-g+1)-2=3n_1-2>0.
   \end{eqnarray*}
\end{enumerate}

\n Now we pass on to prove the first inequality in
\eqref{criteria}. We have
 $$
 C_{12}=n_1n_2(g-1)-n_1d_2+n_2d_1+(d_2-n_2(g-1)-2).
 $$
On the other hand, either $\HH^2_{21}=0$ or else Lemma
\ref{lem:hopf} and Lemma \ref{clif-lem} or Lemma \ref{lem:clif-cs}
imply that
 \begin{equation} \label{eqn:k=3.last}
 \dim \HH^2_{21}\leq h^0(E_1^*\ox K)-1 \leq n_1g -\frac{d_1}{2}
+(n_1-1)\a_i-1.
 \end{equation}
We have the following cases:
 \begin{itemize}
  \item $2n_1-n_2>0$. As we are supposing $n\geq 3$ it follows that $n_1>1$.
Then
   $$
   \qquad \qquad C_{12}=(n_1-1)n_2(g-1) +\a_i(2n_1-n_2)+d_2-2\geq 1+1+1-2
>0.
   $$
   We need to prove that $C_{12}>\dim \HH^2_{21}$ using
   \eqref{eqn:k=3.last}.
   \begin{enumerate}
   \item $n_2=1$, $n_1\geq 2$. Then $C_{12}> n_1g-\frac{d_1}2+(n_1-1)\a_i-1$
    is equivalent to $n_1\a_i +d_2+\frac{d_1}2>n_1+g$.
    In order for coherent systems of type $(1,d_2,2)$ to exist on a Petri
    curve, it is necessary that
    $d_2\geq \frac{g+2}2$. Also
    $$
    d_1>n_1d_2 \geq n_1\frac{g+2}2\geq g+2n_1-2.
    $$
    Easily we get the result.
   \item $n_2=2$, $n_1\geq 2$. Then $C_{12}>n_1g-\frac{d_1}2+(n_1-1)\a_i-1$
    is equivalent to
    $$
    (n_1-2)(g-2)+\a_i(n_1-1)+d_2 +\frac{d_1}2>3.
    $$
    This holds since $d_2\geq 3$, for a stable coherent system of type
    $(2,d_2,2)$ to exist.
   \item $n_2>2$, $n_1\geq 2$. Then generically the map ${\cO}^2\to E_2$ has
no
    kernel, for $(E_2,V_2)\in G_i(n_2,d_2,2)$.
    This happens since in $G_L(n_2,d_2,2)$ all coherent systems have
    this property by Proposition \ref{prop:BGN}, and because all
    $G_i(n_2,d_2,2)$ are birational to each other, by induction hypothesis.
    Therefore the subset $S \subset G_i(n_2,d_2,2)$ of those
    coherent systems $(E_2,V_2)$ such that $\cO^2\to E_2$ is
    not injective is of positive codimension. For $(E_2,V_2)\notin
    S$ we have $\HH^2_{21}=0$ by \eqref{obs-ext2}. So it is enough to prove
    $C_{12}>\dim \HH^2_{21}-1$, i.e.
    $$
    \qquad \qquad\qquad (n_1-2)(g-2)+(n_2-2)(n_1-1)(g-1)+(d_2-\a_i(n_2-2))+
    $$
    $$
    \qquad \qquad \qquad \qquad\qquad \qquad \qquad
+\a_i(n_1-1)+\frac{d_1}2>2.
    $$
    This holds clearly. The only case to be considered
    separately is $g=2$, $n_1=2$, $n_2=3$, $d_1=1$. But in this case
    $\a_i \in \ZZ$, hence $\a_i\geq 1$ and the result follows easily.
   \end{enumerate}

  \item $2n_1-n_2<0$. Then
   $$
   C_{12}=(n_1-1)n_2(g-1) -\a_i(n_2-2n_1)+d_2-2.
   $$
  Now $\a_i=\frac{d_2n_1-d_1n_2}{n_2-2n_1}<\frac{d_2}{n_2-2}$
  gives $\a_i(n_2-2n_1)<d_2-2d_1$. Thus
  $$
  C_{12}>(n_1-1)n_2(g-1)+2d_1-2\geq 0.
  $$
  We have the following cases:
   \begin{enumerate}
   \item $n_1\geq 2$. We use the bound \eqref{eqn:k=3.last}.
    Then  $C_{12}>n_1g -\frac{d_1}2 +(n_1-1)\a_i-1$ is
    equivalent to
    $$
    (n_1-2)(g-2)+(n_2-2)(n_1-1)(g-1)+
    $$
    $$
    \qquad \qquad\qquad \qquad \qquad \qquad
    +(d_2-\a_i(n_2-2n_1))-\a_i(n_1-1)+\frac{d_1}2>3.
    $$
    Use that $\a_i(n_1-1)<d_1$ and $d_2-\a_i(n_2-2n_1)>2d_1$
    to get that the left hand side is bigger or equal than
    $3+2d_1-d_1+\frac{d_1}2>3$.
   \item $n_1=1$, $n_2>2$.
    By Proposition \ref{prop:C21}, $\HH^2_{21}=H^0(E_1^*\ox N_2\ox
    K)^*$, where $N_2\hookrightarrow {\cO}^2\to E_2$. Hence if $N_2=0$
    then $\HH^2_{21}=0$ and we have finished.
    So we may suppose that $N_2\neq 0$. By
    \eqref{eqn:k=3.last} it is enough to prove that
    $C_{12}=d_1n_2-2>g-\frac{d_1}2-1$.
    Let $L$ be the image of ${\cO}^2\to E_2$, which is a line bundle of
    degree $l\geq \frac{g+2}2$. We may write an inclusion
    of coherent systems $(L,V_2) \subset (E_2,V_2)$.
    By $\a_i$-semistability,
    $l+2\a_i\leq \frac{d_2+2\a_i}{n_2}=d_1+\a_i$. So
    $d_1 \geq \a_i+ \frac{g+2}2$ and then
    $$
    d_1n_2-2>n_2\frac{g+2}{2}-2>g-\frac{d_1}2-1.
    $$
   \end{enumerate}
 \end{itemize}
\end{enumerate}
\qed

%%%%%%%%%%%%%%%%%%%%%%%%%%%%%%%%%%%%%%%%%%%%%%%%%%%%%%%%%%%%%%%%%%%
\section{Applications of coherent systems to Brill-Noether
theory}\label{sec:applications}
%%%%%%%%%%%%%%%%%%%%%%%%%%%%%%%%%%%%%%%%%%%%%%%%%%%%%%%%%%%%%%%%%%%

In this section, we shall describe in more detail the relationship
between $G_0(n,d,k)$ and $B(n,d,k)$ introduced in section
\ref{subsec:alpha-small}, and give some applications of our
results to Brill-Noether theory. Although a good deal is known about
non-emptiness of Brill-Noether loci, even quite simple geometrical
properties (for example, irreducibility) have been established only in
a very few cases. The results given here begin to fill these gaps in
our knowledge, and should be regarded as a
sample of what is possible. We plan to return to these questions in
future papers and obtain more extensive and comprehensive results.

Although many of the proofs are valid for all $g$, one may as well
assume in this section that $g\ge2$, since Brill-Noether theory itself is
trivial for $g=0,1$.

%%%%%%%%%%%%%%%%%%%%%%%%%%%%%%%%%%%%%%%%%%%%%%%%%%%%%%%%%%%%%%%%%%%%%%%%
\subsection{General remarks.}\label{subsec:bngen}
%%%%%%%%%%%%%%%%%%%%%%%%%%%%%%%%%%%%%%%%%%%%%%%%%%%%%%%%%%%%%%%%%%%%%%%
\begin{lemma}\label{lem:bn1}
If $\beta(n,d,k)\ge n^2(g-1)+1$, then $B(n,d,k)=M(n,d)$.
\end{lemma}

\noindent{\em Proof. \/} By (\ref{bnnumber}),
\begin{equation}\label{eq:bn1}
\beta(n,d,k)\ge n^2(g-1)+1\Leftrightarrow d-n(g-1)\ge k.
\end{equation}
When these equivalent conditions hold, it follows from the
Riemann-Roch Theorem that, for any $E\in M(n,d)$,
 $$
 \dim H^0(E)\ge k.
 $$
So $B(n,d,k)=M(n,d)$. \hfill $\Box$

\medskip On the other hand, we have
\begin{lemma}\label{lem:bn2}
If $\beta(n,d,k)\le n^2(g-1)$, then every irreducible component $B$
of $B(n,d,k)$ contains a point outside $B(n,d,k+1)$.
\end{lemma}

\noindent{\em Proof. \/} (This is \cite[Lemma 2.6]{Lau}; for the
convenience of the reader, we include a proof.) The content of the
statement is that there exists $E\in B$ such that $\dim H^0(E)=k$.
To see this, note first that, if $\dim H^0(E')\ge1$ and $P$ is a
point of $X$ such that the sections of $E'$ generate a non-zero
subspace of the fibre $E'_P$, we can find an extension
 $$
 0\to F\to E'\to{\mathcal O}_P\to 0
 $$
such that the map $H^0(E')\to {\mathcal O}_P$ is non-zero and
hence $\dim H^0(F)=\dim H^0(E')-1$.

Now let $E'$ be a point of $B$ not contained in any other
irreducible component of $B(n,d,k)$ and suppose that $H^0(E')=k+r$
with $r\ge1$. By iterating the above construction, we can find
points $P_1,\ldots,P_r$ of $X$ and an exact sequence
 $$
 0\to F\to E'\to{\mathcal O}_{P_1}\oplus\ldots\oplus{\mathcal
 O}_{P_r}\to 0
 $$
such that $\dim H^0(F)=\dim H^0(E')-r=k$. Now consider the extensions
 \begin{equation}\label{eq:bn2}
 0\to F\to E\to{\mathcal O}_{Q_1}\oplus\ldots\oplus{\mathcal
 O}_{Q_r}\to 0,
 \end{equation}
where $Q_1,\ldots,Q_r\in X$.
These form an irreducible family of bundles with $\dim H^0(E)\ge k$,
whose generic member is stable (since $E'$ is stable). It follows
that the generic extension (\ref{eq:bn2}) belongs to $B$. Moreover,
by the Riemann-Roch Theorem and (\ref{eq:bn2}),
 \begin{eqnarray*}
 \dim H^1(F)&=&\dim H^0(F)-(d-r)+n(g-1)\\
 &>& k-k+r=r.
 \end{eqnarray*}
By considering the dual sequence
 $$
 0\to E^*\otimes K\to F^*\otimes K\to{\mathcal
 O}_{Q_1}\oplus\ldots\oplus{\mathcal O}_{Q_r}\to 0,
 $$
in which $\dim H^0(F^*\otimes K)>r$, we can choose
$Q_1,\ldots,Q_r$ and $E$ so that
 $$
 \dim H^0(E^*\otimes K)=\dim H^0(F^*\otimes K)-r;
 $$
hence (again by Riemann-Roch)
 $$
 \dim H^0(E)=\dim H^0(F).
 $$
It follows that the generic extension (\ref{eq:bn2}) satisfies $\dim
H^0(E)=k$. Since we already know that $E\in B$, this completes the
proof.\hfill $\Box$

\medskip As envisaged at the end of section \ref{subsec:alpha-small}, we
introduce
\begin{cond}\label{cond} {\ }
\begin{itemize}
\item $\beta(n,d,k)\le n^2(g-1)$,
\item $G_0(n,d,k)$ is irreducible,
\item $B(n,d,k)\ne\emptyset$.
\end{itemize}
\end{cond}

For the moment we do not assume that ${\GCD}(n,d,k)=1$ or that
$G_0(n,d,k)$ is smooth. We denote by
 $$
 \psi:G_0(n,d,k)\to\widetilde{B}(n,d,k)
 $$
the map given by assigning to every $(E,V)\in G_0(n,d,k)$ the
underlying bundle $E$ (see (\ref{bnmap})).

\begin{theorem}\label{thm:bn1}
Suppose Conditions \ref{cond} hold. Then
\begin{enumerate}
\item $B(n,d,k)$ is irreducible,
\item $\psi$ is one-to-one over $B(n,d,k)-B(n,d,k+1)$,
\item $\dim B(n,d,k)=\dim G_0(n,d,k)$,
\item for any $E\in B(n,d,k)-B(n,d,k+1)$, the linear map
 $$
 {\rm d}{\psi}:T_{(E,H^0(E))}G_0(n,d,k)\longrightarrow
 T_EB(n,d,k)
 $$
of Zariski tangent spaces is an isomorphism.\end{enumerate}
\end{theorem}
\n{\em Proof. \/} (i) If $E\in B(n,d,k)$, then $(E,V)\in G_0(n,d,k)$ for
any $k$-dimensional subspace of $H^0(E)$. It follows that the image
of $\psi$ contains $B(n,d,k)$ as a non-empty Zariski-open subset.
Since $G_0(n,d,k)$ is irreducible, it follows that $B(n,d,k)$ is
irreducible.

(ii) If $E\in B(n,d,k)-B(n,d,k+1)$, then
$\psi^{-1}(E)=\{(E,H^0(E)\}$.

(iii) follows from (i), (ii) and Lemma \ref{lem:bn2}.

(iv) Taking $(E',V')=(E,V)$ in (\ref{long-exact}) and putting $V=H^0(E)$,
we get a map
 $$
 {\Ext}^1((E,H^0(V)),(E,H^0(V)))\to {\Ext}^1(E,E)
 $$
which can be identified with the map
 $$
 T_{(E,H^0(E))}G_0(n,d,k)\to T_EM(n,d)
 $$
induced by $\psi$. By (\ref{long-exact}) this map is injective and
its image is
 $$
 {\Ker}({\Ext}^1(E,E)\to {\Hom}(H^0(E),H^1(E))).
 $$
By standard Brill-Noether theory, this image becomes identified with
the subspace $T_EB(n,d,k)$ of $T_EM(n,d)$.\hfill$\Box$

\begin{corollary}\label{cor:smooth} Suppose Conditions \ref{cond} hold
and $G_0(n,d,k)$ is smooth. Then $\psi$ is an isomorphism over
$B(n,d,k)-B(n,d,k+1)$. Moreover, if ${\GCD}(n,d,k)=1$, then
$G_0(n,d,k)$ is a desingularisation of the closure
$\overline{B(n,d,k)}$ of $B(n,d,k)$ in the projective variety
$\widetilde{M}(n,d)$.
\end{corollary}

\n{\em Proof. \/} The first part follows from (ii) and (iv). For
the second part, recall that, when ${\GCD}(n,d,k)=1$, $G_0(n,d,k)$
is projective; hence the image of $\psi$ is precisely
$\overline{B(n,d,k)}$.\hfill $\Box$

\begin{corollary}\label{cor:coprime} Suppose Conditions \ref{cond} hold,
$G_0(n,d,k)$ is smooth and $(n,d)=1$. Then $B(n,d,k)$ is projective
and $G_0(n,d,k)$ is a desingularisation of $B(n,d,k)$.
\end{corollary}

\n{\em Proof. \/} In this case $M(n,d)=\widetilde{M}(n,d)$.\hfill $\Box$

%%%%%%%%%%%%%%%%%%%%%%%%%%%%%%%%%%%%%%%%%%%%%%%%%%%%%%%%%%%%%%%%%
\subsection{Irreducibility and dimension of
Brill-Noether loci.}\label{sec:irred}
%%%%%%%%%%%%%%%%%%%%%%%%%%%%%%%%%%%%%%%%%%%%%%%%%%%%%%%%%%%%%%%%%

In many cases our methods yield information about the irreducibility
and dimension of $B(n,d,k)$, and more precisely about its birational
structure. We illustrate this with results for $k=1,2,3$, where we
have good estimates for the codimensions of the flips. The main
respect in which our results improve those previously known is that
they impose no restriction on $d$ other than that required for
the Brill-Noether locus to be non-empty and not equal to $M(n,d)$.

We begin with $k=1$.

\begin{theorem}\label{thm:bn2}
Suppose $0<d\le n(g-1)$. Then
\begin{enumerate}
\item $G_0(n,d,1)$ is a desingularisation of
$\overline{B(n,d,1)}$,
\item $B(n,d,1)$ is irreducible of dimension
$\beta(n,d,1)$, smooth outside $B(n,d,2)$,
\item $B(n,d,1)$ is birationally equivalent to a fibration
over $M(n-1,d)$ with fibre ${\mathbb P}^{d+(n-1)(g-1)-1}$,
\item if $(n-1,d)=1$, $B(n,d,1)$ is birationally equivalent to
$$M(n-1,d)\times{\mathbb P}^{d+(n-1)(g-1)-1}.$$
\end{enumerate}
\end{theorem}

\begin{remark}\begin{em}\label{bn:rem1}
In the case $d=n(g-1)$, a stronger form of (i) is proved in \cite{RV}.
Part (ii) is proved in \cite{Su}. Parts (iii) and (iv) are implicit in
\cite{Su}. We have chosen to prove the complete theorem to illustrate
our methods.
\end{em}\end{remark}

\n{\em Proof. \/} We first check Conditions \ref{cond}. The first
follows at once from (\ref{eq:bn1}), the second from Theorem
\ref{thm:k=1} and the third is elementary and well known (see for
example \cite{Su}). Moreover $G_0(n,d,1)$ is smooth of dimension
$\beta(n,d,1)$ by Theorem \ref{thm:k=1} (or Proposition
\ref{prop:smooth}).

Parts (i) and (ii) now follow from Theorem \ref{thm:bn1} and
Corollary \ref{cor:smooth}. Part (iii) follows from Theorem
\ref{G_L(k<n)} and Theorem \ref{thm:k=1}, as does part (iv) if we
note that in this case the existence of a universal bundle over
$X\times M(n-1,d)$ implies that the fibration of Theorem
\ref{G_L(k<n)} is locally trivial in the Zariski topology.\hfill
$\Box$

For $k=2,3$, we need a lemma.

\begin{lemma}\label{lem:bn3} Suppose $k=2$ or $3$, $n\ge2$ and that
$X$ is a Petri curve of genus $g\ge2$. Then $B(n,d,k)$ is non-empty
precisely in the following cases:
\renewcommand{\theenumi}{\alph{enumi}}
\begin{enumerate}
\item $k=2$, $n=2$, $d\ge3$,
\item $k=2$, $n\ge3$, $d\ge1$,
\item $k=3$, $n=2$, $d\ge\frac{2g+6}{3}$,
\item $k=3$, $n=3$, $d\ge4$,
\item $k=3$, $n=4$, either $g=2$ and $d\ge2$ or $g\ge3$ and
$d\ge1$,
\item $k=3$, $n\ge5$, $d\ge1$.
\end{enumerate}
\renewcommand{\theenumi}{\roman{enumi}}
\end{lemma}

\begin{remark}\begin{em} The Petri condition is required only for case
(c).
\end{em}\end{remark}

\n{\em Proof. \/} All parts except (c) follow from \cite{BGN} and either
\cite{T1} or \cite{M5}. For (c), see \cite{Bu2} or \cite{Tan}.\hfill $\Box$

\begin{theorem}\label{thm:bn3} Let $X$ be a Petri curve of genus $g\ge2$,
$k=2$ or $3$, $d<n(g-1)+k$. Suppose further that one of the
conditions of Lemma \ref{lem:bn3} holds. Then

{\rm(i)} $B(n,d,k)$ is irreducible of dimension $\beta(n,d,k)$.

\n If in addition $k<n$ (i.e. in cases {\em (b), (e), (f)} of Lemma
\ref{lem:bn3}), then

{\rm(ii)}  $B(n,d,k)$ is birationally equivalent to a fibration
over $M(n-k,d)$ with fibre ${\Gr}(k,d+(n-k)(g-1))$;

{\rm(iii)} if $(n-k,d)=1$, $B(n,d,k)$ is birationally equivalent
to  $$M(n-k,d)\times{\Gr}(k,d+(n-k)(g-1)).$$
\end{theorem}

\n{\em Proof. \/} (i) Note that the conditions of Lemma \ref{lem:bn3} for
the non-emptiness of $B(n,d,k)$ are exactly the same as those of
Theorems \ref{thm:k=2} and \ref{thm:k=3} for the non-emptiness of
$G_0(n,d,k)$.
Conditions
\ref{cond} follow from (\ref{eq:bn1}) and Theorems \ref{thm:k=2} and
\ref{thm:k=3}. The
result now follows from Theorem \ref{thm:bn1}.

(ii) and (iii) follow from Theorem \ref{G_L(k<n)} in the same way as the
corresponding parts of Theorem \ref{thm:bn2}.\hfill $\Box$

\begin{remark}\begin{em} When $k<n$ (cases (b), (e), (f)), the
irreducibility of $B(n,d,k)$ for $d<\min\{2n,n+g\}$
and the fact that $B(n,d,k)$ has the expected dimension
for $d\le 2n$ have been proved previously \cite{BGN, M1}.
For $k\le n$ (all cases except (c)),
it was proved in \cite{T1} that $B(n,d,k)$ has a component of the
expected dimension. Parts (ii) and (iii) are known for $d<\min\{2n,n+g\}$
\cite{M1}.
\end{em}\end{remark}

%%%%%%%%%%%%%%%%%%%%%%%%%%%%%%%%%%%%%%%%%%%%%%%%%%%%%%%%%%%%%%%%%%%
\subsection{Picard group.}\label{sec:pic}
%%%%%%%%%%%%%%%%%%%%%%%%%%%%%%%%%%%%%%%%%%%%%%%%%%%%%%%%%%%%%%%%%%
Our methods become potentially even more useful in computing
cohomological information about Brill-Noether loci. In general,
the calculations will be complicated and we restrict attention here
to  computing the Picard group in the case $k=1$.

\begin{theorem}\label{thm:bn4} Let $X$ be a Petri curve of genus $g\ge2$.
Suppose
$0<d\le n(g-1)$, $n\ge3$, $(n-1,d)=1$ and $(n,d)=1$. Then
 $$
 {\Pic}(B(n,d,1)-B(n,d,2))\cong{\Pic}(M(n-1,d))\times{\mathbb Z}.
 $$
\end{theorem}

\n{\em Proof. \/} Note
first that, by Theorem \ref{G_L(k<n)}, $G_L(n,d,1)$ is a projective bundle
over
$M(n-1,d)$, so
 $$
 {\Pic}(G_L(n,d,1))={\Pic}(M(n-1,d))\times{\mathbb Z}.
 $$
{}From the proof of Theorem \ref{thm:k=1}, we see that the
codimensions $C_{12}$, $C_{21}$ are both at least $2$ (we need
$n\ge3$ here since otherwise we could have $n_1=n_2=1$, $d_2=1$,
giving $C_{21}=1$). Hence
 $$
 {\Pic}(G_0(n,d,1))={\Pic}(M(n-1,d))\times{\mathbb Z}.
 $$

To complete the proof, we need to show that $\psi^{-1}(B(n,d,2))$
has codimension at least $2$ in $G_0(n,d,1)$. Now the fibre of
$\psi$ over a point of $B(n,d,k)-B(n,d,k+1)$ is a projective space of
dimension $k-1$. It is therefore sufficient to prove that $B(n,d,k)$
has codimension at least $k+1$ in $B(n,d,1)$ for all $k\ge2$. In view
of Lemma \ref{lem:bn2}, it is enough to prove this for $k=2$, i.e.
to prove
 $$
 {\codim}_{B(n,d,1)}B(n,d,2)\ge3.
 $$

For this, note that
\begin{eqnarray*}
\beta(n,d,2)&=&\beta(n,d,1)-n(g-1)+d-3\\
&\le&\beta(n,d,1)-3
\end{eqnarray*}
since $d\le n(g-1)$. The result now follows from Theorem
\ref{thm:bn3}.

\begin{remark}\begin{em}\label{rem:bn2}
Note that we need Theorem \ref{thm:bn3} here to show that
$B(n,d,2)$ always has the expected dimension. This is the only
point in the proof where the Petri condition is used. It may be
that this condition is not essential.
\end{em}\end{remark}

%%%%%%%%%%%%%%%%%%%%%%%%%%%%%%%%%%%%%%%%%%%%%%%%%%%%%%%%%%%%%%%%%%

\end{document}